\documentclass[leqno]{macrorend}

\volumeyear{xxxx}\yearnumber{x}\volumenumber{xx}

\usepackage[latin1]{inputenc}
\usepackage[PostScript]{diagrams}
\diagramstyle[labelstyle=\scriptstyle]
\def\CD{\diagram[amstex]}
\newcommand{\wA}{\widetilde A}

\newcommand{\pr}{\operatorname{pr}}
\newcommand{\inj}{\operatorname{i}}

\newcommand{\rnp}{{\mathbb R}^n_+}

\newcommand{\crnp}{\overline{\mathbb R}^n_+}

\newcommand{\comega}{\overline\Omega }

\newcommand{\ang}[1]{\langle {#1} \rangle}

\newcommand{\simto}{\overset\sim\rightarrow}

\newcommand{\Ami}{A_{\min}}
\newcommand{\Ama}{A_{\max}}

\begin{document}

% select a language among: english, french, italian
%
\selectlanguage{english}

% if you dont have footnote, cancel \footnotemark[1]
% separate several authors with ' - '
%
\articolo[Krein resolvent formulas]{Krein resolvent formulas for 
 elliptic
 boundary problems in nonsmooth domains}{G.~Grubb}%\footnotemark[1]}

% text of the note (number from 1 to 9)
%
%\footnotetext[1]{Text of the footnote.}
%%gik ikke\dedicatory{Dedicated to Luigi Rodino on the occasion 
%%of his 60-th birthday.}
\begin{abstract}
The paper reports on a recent construction of $M$-functions and Kre\u\i{}n
resolvent formulas for general closed extensions of an adjoint pair, and
their implementation to boundary value problems for second-order strongly elliptic
operators on smooth domains.
The results are then extended to domains with $C^{1,1}$ H\"older smoothness,
by use of a recently developed calculus of pseudodifferential boundary operators with
nonsmooth symbols.
\end{abstract}

\section{Introduction.}\label{Section0}
In the study of boundary value problems for ordinary differential
equations, the Weyl-Titchmarsh $m$-function has played an important
role for many years; it allows a reduction of questions concerning the
resolvent $(\wA-\lambda )^{-1}$ of a realisation $\wA$ to questions
concerning an associated family $M(\lambda )$ of matrices,
holomorphic in $\lambda \in \varrho (\wA)$. Moreover, there is a formula describing the
difference between %\linebreak
$(\wA-\lambda )^{-1}$ and the resolvent of a well-known
reference problem in terms of $M(\lambda )$, a so-called
Kre\u\i{}n resolvent formula. The concepts have also been introduced
in connection with the abstract theories of extensions of symmetric
operators or adjoint pairs in Hilbert spaces, initiated by Kre\u\i{}n
\cite{K47} and  Vishik \cite{V52}. The literature on this is
abundant, and we refer to e.g.\ Brown, Marletta, Naboko and Wood
\cite{BMNW08} and Brown, Grubb and Wood \cite{BGW08} for accounts
of the development, and references. For elliptic partial differential
equations  in higher
dimensions, concrete interpretations of
$M(\lambda )$ have been taken up
in recent years, e.g.\ in Amrein and Pearson
\cite{AP04}, Behrndt and Langer \cite{BL07}, and in \cite{BMNW08}; here $M(\lambda )$ is a
family of operators defined over the boundary. 
In the present paper we report on the
latest development in nonsymmetric cases worked out in
\cite{BGW08}; it uses the early work of Grubb \cite{G68} as an
important ingredient. 

The interest of this in a context of pseudodifferential
operators  is that $M(\lambda )$ in elliptic cases, 
and also in some
nonelliptic cases, is a pseudodifferential operator ($\psi $do), to
which $\psi$do  methods can be applied. The new results in the present
paper are concerned with situations with a nonsmooth boundary. Our strategy
here is to apply the nonsmooth pseudodifferential boundary operator 
($\psi $dbo) calculus introduced by Abels
\cite{A05}. We show that when the domain is $C^{1,1}$ and the
given strongly elliptic second-order operator $A$ has smooth coefficients, then indeed the
$M$-function can be defined as a generalized $\psi$do over the
boundary, and a Kre\u\i{}n formula holds.
Selfadjoint cases have been treated
under various nonsmoothness hypotheses in Gesztesy and Mitrea
\cite{GM08}, Posilicano and Raimondi \cite{PR08}, but
the present study allows nonselfadjoint operators, and includes a
discussion of Neumann-type boundary conditions. Besides bounded
domains, we also treat
exterior domains and perturbed
halfspaces. 

The author thanks Helmut Abels for useful conversations.

%Y and then we establish 
%$some results on the implementation to boundary value 
%problems for second-order 
%strongly elliptic considered on a

\section{Abstract results.}\label{Section1}

We begin by recalling the theory of extensions and $M$-functions
established in works of Brown, Wood and the author \cite{BGW08} and \cite{G68}.

There is given an adjoint pair of closed, densely defined linear
operators $A_{\min}$, $A'_{\min}$ in a Hilbert
space $H$:
\[
A_{\min}\subset (A'_{\min})^*=A_{\max},\quad A'_{\min}\subset (A_{\min})^*=A'_{\max}.
\]
Let $\mathcal M$ denote the set of linear operators lying between the
minimal and maximal operator:
\[\mathcal M=\{\wA\mid \Ami\subset
\wA\subset \Ama\},\quad
\mathcal M'=\{\wA'\mid \Ami'\subset
\wA'\subset \Ama'\}.\]
 Write $\wA u$ as $Au$ for any $\wA$, and $\wA' u$ as $A'u$ for any $\wA'$.
Assume that there exists an $A_\gamma \in\mathcal M$ with $0\in
\varrho (A_\gamma )$; then $A_\gamma ^*\in \mathcal M'$ with $0\in \varrho
(A_\gamma ^*)$. 
We shall define $M$-functions
for {\it any} closed $\wA\in
\mathcal M$. 

First recall some details from the treatment of extensions in
\cite{G68}: Denote
\[
Z=\operatorname{ker}\Ama, \quad Z'=\operatorname{ker}\Ama'.
\]
Define the basic non-orthogonal decompositions 
\begin{align}
D(\Ama)&=D(A_\gamma )\dot+ Z,\text{ denoted }u=u_\gamma +u_\zeta
=\pr_\gamma u+\pr_\zeta u,\nonumber\\
D(\Ama')&=D(A_\gamma ^*)\dot+ Z',\text{ denoted } v=v_{\gamma '}+v_{\zeta '}
=\pr_{\gamma '}v+\pr_{\zeta '}v;\nonumber
\end{align}
here $\pr_\gamma =A_\gamma ^{-1}\Ama$, $\pr_\zeta =I-\pr_\gamma $, and
$\pr_{\gamma '}=(A^*_\gamma )^{-1}\Ama'$, $\pr_{\zeta '}=I-\pr_{\gamma '}$.
By $\pr_Vu=u_V$ we denote the {\it orthogonal projection} of $u$ onto
$V$.

The following ``abstract Green's formula'' holds:
\begin{equation}(A u,v) - (u, A' v)=((A u)_{Z'}, v_{\zeta '})-(u_\zeta ,
(A 'v)_Z).\label{tag1.1}
\end{equation}
It can be used to show that when $\wA\in \mathcal M$ and we set
$W=\overline{\pr_{\zeta '} D(\widetilde A^*)}$, then
\[
\{\{u_\zeta , (A u)_W\}\mid u\in D(\wA)\}\text{ is a graph.} 
\]
Denoting
the operator with this graph by $T$, we have:

\begin{theorem}\label{Theorem1.1} {\rm {\cite{G68}}} For the closed $\wA\in\mathcal M$,
there is a {\rm 1--1} correspondence
\[
\wA \text{ closed } \longleftrightarrow \begin{cases} T:V\to W,\text{ closed, densely
defined}\\
\text{with }V\subset Z,\; W\subset Z',\text{ closed subspaces.} \end{cases}
\]
\noindent Here $D(T)=\pr_\zeta  D(\widetilde A)$, $V=\overline{D(T)}$, $
W=\overline{\pr_{\zeta '} D(\widetilde A^*)}$, and
\begin{equation}
Tu_\zeta=(Au)_W \text{ for all }u\in D(\wA),\text{ (the {\bf defining
equation})}.\label{tag1.2}
\end{equation}
In this correspondence,

{\rm (i)} $\wA^*$ corresponds similarly to $T^*:W\to V$.

{\rm (ii)} $\operatorname{ker}\wA=\operatorname{ker}T$;
\quad $\operatorname{ran}\wA=\operatorname{ran}T+(H\ominus W)$.

{\rm (iiii)} When $\wA$ is invertible, 
\[\wA^{-1}=A_\gamma
^{-1}+\inj_{V\to H}T^{-1}\pr_W.\]

\end{theorem}

Here $\inj_{V\to H}$ indicates the injection of $V$ into $H$ (it is
often left out). 

Now provide the operators with a spectral parameter $\lambda $, then this
implies, with
\begin{align}
Z_\lambda& =\operatorname{ker}(\Ama-\lambda ), \quad Z'_{\bar\lambda }=\operatorname{ker}(\Ama'-\bar\lambda ),\nonumber\\
D(\Ama)&=D(A_\gamma )\dot+ Z_\lambda ,  \quad u=u^\lambda _\gamma +u^\lambda _\zeta
=\pr^\lambda _\gamma u+\pr^\lambda _\zeta u,\text{ etc.:}\nonumber
\end{align} 

\begin{corollary}\label{Corollary1.2} Let $\lambda \in \varrho (A_\gamma )$. For the
closed $\wA\in\mathcal M$,  there is a {\rm 1--1} correspondence
\[
\wA -\lambda  \longleftrightarrow \begin{cases} T^\lambda :V_\lambda \to W_{\bar\lambda },\text{ closed, densely
defined}\\
\text{with }V_\lambda \subset Z_\lambda ,\; W_{\bar\lambda }\subset Z'_{\bar\lambda },\text{ closed subspaces.} \end{cases}
\]
\noindent Here $D(T^\lambda )=\pr^\lambda _\zeta  D(\widetilde A)$, $V_\lambda =\overline{D(T^\lambda )}$, $
W_{\bar\lambda }=\overline{\pr^{\bar\lambda }_{\zeta '} D(\widetilde A^*)}$, and
\[
T^\lambda u^\lambda _\zeta=((A-\lambda )u)_{W_{\bar\lambda }} \text{ for all }u\in D(\wA).\]

 Moreover,

{\rm (i)} $\operatorname{ker}(\wA-\lambda )=\operatorname{ker}T^\lambda $;
\quad $\operatorname{ran}(\wA-\lambda )=\operatorname{ran}T^\lambda +(H\ominus W_{\bar\lambda })$.

{\rm (ii)} When $\lambda \in \varrho (\wA)\cap \varrho (A_\gamma )$, 
\begin{equation}(\wA-\lambda
)^{-1}=(A_\gamma -\lambda )
^{-1}+\inj_{V_\lambda \to H}(T^\lambda )^{-1}\pr_{W_{\bar\lambda
  }}.\label{tag1.3}
\end{equation}

\end{corollary}

This gives a Kre\u\i{}n resolvent formula for any closed
$\wA\in\mathcal M$. 

The operators $T$ and $T^\lambda $ are related in the following way: Define
\begin{align}
E^\lambda &=I+\lambda (A_\gamma -\lambda )^{-1},\quad
F^\lambda =I-\lambda A_\gamma ^{-1},\nonumber\\
E^{\prime\bar\lambda }&=I+\bar\lambda (A^*_\gamma -\bar\lambda )^{-1},\quad
F^{\prime\bar\lambda }=I-\bar\lambda (A^*_\gamma )^{-1},\nonumber
\end{align}
then $E^\lambda F^\lambda =F^\lambda E^\lambda =I$,
$E^{\prime\bar\lambda }F^{\prime\bar\lambda }=F^{\prime\bar\lambda
}E^{\prime\bar\lambda }=I$ on $H$. Moreover, $E^\lambda $ and
$E^{\prime\bar\lambda }$ restrict to 
homeomorphisms 
\[
E^\lambda _V: V\simto V_\lambda ,\quad E^{\prime\bar\lambda
}_W:W\simto W_{\bar\lambda },\]
with inverses denoted $F^\lambda _V$ resp.\ $F^{\prime\bar\lambda }_W$.
In particular, $D(T^\lambda )=E^\lambda _V D(T)$.

\begin{theorem}\label{Theorem1.3}
 Let $
G^\lambda _{V,W}=-\pr_W\lambda
E^\lambda \inj_{V\to H}$; then 
\begin{equation}(E^{\prime\bar\lambda }_W)^* T^\lambda E^\lambda _V=T+G^\lambda
_{V,W}.\label{tag1.2}
\end{equation}
In other words, $T$ and $T^\lambda $ are related by the
commutative diagram (where the horizontal maps are homeomorphisms)
\[
\CD
V_\lambda     @<{E^\lambda _{V}} <<  V    \\
@V  T^\lambda VV           @VV  T+G^\lambda _{V,W} V\\
W_{\bar\lambda }   @>>(E^{\prime\bar\lambda }_W)^* >   W
\endCD \hskip2cm D(T^\lambda )=E^\lambda _V D(T).
 \]
\end{theorem}
 
This is a straightforward elaboration of \cite{G74}, Prop.\ 2.6.

Now let us introduce boundary triplets and $M$-functions.
The general setting is the following: There is given a pair of 
Hilbert spaces $\mathcal H$, $\mathcal K$ and two
pairs of
``boundary operators''
\begin{equation}
 \begin{pmatrix}\Gamma_1\\ \quad\\\Gamma_0 \end{pmatrix}: D(\Ama)\to  \begin{matrix}{\mathcal H}\\\times\\{\mathcal K}\end{matrix}, \quad
\begin{pmatrix}\Gamma'_1\\ \quad\\ \Gamma'_0 \end{pmatrix}: D(\Ama')\to
\begin{matrix}{\mathcal K}\\ \times\\{\mathcal H}\end{matrix},
\label{tag1.3}
\end{equation}
bounded with respect to the graph norm and  surjective, such that
\[
D(\Ami)=D(\Ama)\cap \operatorname{ker}\Gamma _1\cap
\operatorname{ker}\Gamma _0,\quad 
D(\Ami')=D(\Ama')\cap \operatorname{ker}\Gamma '_1\cap
\operatorname{ker}\Gamma '_0,
\] and  for all $u\in D(\Ama)$, $v\in D(\Ama')$,
\[
(A u,v) - (u, A' v)=(\Gamma _1u, \Gamma '_0v)_{\mathcal H}-(\Gamma _0u ,
\Gamma '_1v)_{\mathcal K}.
\]
Then the three pairs $\{\mathcal H, \mathcal K\}$, $\{\Gamma _1,\Gamma
_0\}$ and $\{\Gamma '_1,\Gamma '
_0\}$ are said to form a {\it boundary triplet}. (See \cite{BMNW08} and
\cite{BGW08} for references to the literature on this.)

Note that under our assumptions,
the choice
\begin{equation}
\mathcal H=Z',\quad \mathcal K=Z,\quad\begin{pmatrix} \Gamma _1u\\
  \Gamma _0u
\end{pmatrix}=\begin{pmatrix} (A u)_{Z'}\\
u_\zeta \end{pmatrix},\quad
\begin{pmatrix} \Gamma '_1v\\ \Gamma '_0v\end{pmatrix} =\begin{pmatrix} (A' v)_{Z}\\
v_{\zeta '}\end{pmatrix},\label{tag1.3a}
\end{equation}
defines a boundary triplet, cf.\ \eqref{tag1.1}.

Following \cite{BMNW08}, the boundary triplet is used to define 
operators $A_T\in \mathcal M$ and $A'_{T'}\in\mathcal M'$ for any pair of
operators $T\in \mathcal L(\mathcal K,\mathcal H)$, $T'\in \mathcal L(\mathcal H,\mathcal K)$
by
\begin{equation}
D(A_T)=\operatorname{ker}(\Gamma _1-T\Gamma _0),\quad
D(A'_{T'})=\operatorname{ker}(\Gamma '_1-T'\Gamma '_0).\label{tag1.3b}
\end{equation}
Then they show:

\begin{proposition}\label{Proposition1.3a} For $\lambda \in \varrho (A_T)$, there is a well-defined
$M$-function $M_T(\lambda )$ determined  by
\[
M_T(\lambda ): \operatorname{ran}(\Gamma _1-T\Gamma _0)\to {\mathcal
K},\quad M_T(\Gamma _1-T\Gamma _0)u=\Gamma _0u\text{ for all }u\in Z_\lambda . 
\]
Likewise, for $\lambda \in \varrho (A'_{T'})$, the
function $M'_{T'}(\lambda )$ is determined similarly by
\[
M'_{T'}(\lambda ): \operatorname{ran}(\Gamma '_1-T'\Gamma '_0)\to {\mathcal
H},\quad M'_{T'}(\Gamma '_1-T'\Gamma '_0)v=\Gamma '_0v\text{ for all }v\in Z'_\lambda . 
\]
Here, when $\varrho (A_T)\ne \emptyset$,
\[
(A_T)^*=A'_{T^*}.
\]
\end{proposition}

This was set in relation to Theorem \ref{Theorem1.1} in \cite{BGW08}:  Take the
boundary triplet defined in \eqref{tag1.3a}. Then the formula for $D(A_T)$ in
\eqref{tag1.3b} is the same as the defining equation \eqref{tag1.2} for $D(\wA)$. 
For the sake of generality,
 allow also unbounded, densely defined, closed operators $T:Z\to
Z'$; then in fact the formulas in Proposition \ref{Proposition1.3a} still lead to a well-defined
$M$-function $M_{T}(\lambda )$. We denote $A_T$ by $\wA$ and
$M_T(\lambda )$ by $M_{\wA}(\lambda )$,
when they come from the special choice \eqref{tag1.3a} of boundary triplet. Then
we have: 
%explained in \cite{BMNW08}, 

\begin{theorem}\label{Theorem1.4}  
Let $\wA$ correspond to $T:Z\to Z'$ by
Theorem {\rm \ref{Theorem1.1}}. For any $\lambda \in \varrho (\wA)$,
$M_{\wA}(\lambda )$ is in $\mathcal L(Z',Z)$ and satisfies
\[
M_{\wA}(\lambda )=\pr_\zeta (I-(\wA-\lambda )^{-1}(\Ama-\lambda ))A_\gamma
^{-1}\inj_{Z'\to H}.\]
\noindent {\bf Moreover,} $M_{\wA}(\lambda )$ relates to $T$ and $T^\lambda $
by:
\begin{equation}M_{\wA}(\lambda )=-(T+G^\lambda _{Z,Z'})^{-1}=
- F_Z^\lambda(T^\lambda )^{-1}(F_{Z'}^{\prime\bar\lambda })^*\text{, for }\lambda \in \varrho (\wA)\cap
\varrho (A_\gamma ).\label{tag1.4}
\end{equation}
\end{theorem}

This takes care of those operators $\wA$ for which $\pr_\zeta
D(\widetilde A)$ is dense in $Z$ and $\pr_{\zeta '} D(\widetilde A^*)$
is dense in $Z'$. But the construction extends in a natural way to all the closed
$\wA\in\mathcal M$, giving the following result:

\begin{theorem}\label{Theorem1.5}   Let $\wA$ correspond to $T:V\to W$ by Theorem
{\rm \ref{Theorem1.1}}. For any
$\lambda \in \varrho (\wA)$, there is a well-defined
$M_{\wA}(\lambda )\in \mathcal L(W,V)$, holomorphic in $\lambda $
and satisfying

{\rm (i)} $M_{\wA}(\lambda )=\pr_\zeta (I-(\wA-\lambda )^{-1}(\Ama-\lambda ))A_\gamma
^{-1}\inj_{W\to H}.$

{\rm (ii)} When $\lambda \in \varrho (\wA)\cap
\varrho (A_\gamma )$, 
\[
M_{\wA}(\lambda )=-(T+G^\lambda _{V,W})^{-1}.
\]

{\rm (iii)} For $\lambda \in \varrho (\wA)\cap
\varrho (A_\gamma )$, it enters in a Kre\u\i{}n resolvent formula  
\begin{equation}
(\wA-\lambda
)^{-1}=(A_\gamma -\lambda )
^{-1}-\inj_{V_\lambda \to H}E^\lambda _VM_{\wA}(\lambda
)(E^{\prime\bar\lambda }_W)^*\pr_{W_{\bar\lambda }}.\label{tag1.5a}
\end{equation}

\end{theorem}

Other Kre\u{\i}n-type resolvent formulas in a general framework of
{\it relations} can be found in Malamud and Mogilevski\u\i{} 
\cite[Section 5.2]{MM02}.

\section{Neumann-type conditions for
second-order operators.}\label{Section2}

The abstract theory can be applied to elliptic realisations
by use of suitable mappings going to and from the boundary, allowing
an interpretation in terms of boundary conditions. We shall demonstrate
this in the strongly elliptic second-order case.

Let $\Omega $ be an open subset of ${\Bbb R}^n$ of one of the
following three types: 1) $\Omega $ is bounded, 2) $\Omega $ 
is the
complement of a bounded set (i.e., is an
exterior domain), 
 or 3) there is a ball $B(0,R)$ with
center $0$ and radius $R$ such that $\Omega \setminus
B(0,R)=\rnp\setminus B(0,R)$ (we then call $\Omega$ a perturbed
halfspace).
More general sets or manifolds could be considered in a similar way, namely
the so-called admissible manifolds as defined in the book \cite{G96}. 

The sets will in the present section be assumed to be $C^\infty $;
later from Section \ref{Section4} on they will be taken to be 
$C^{k,\sigma }$, where
$k$ is an integer $\ge 0$ and $ \sigma \in \,]0,1]$. (Recall that the
norm on the H\"older space 
$C^{k,\sigma }(V)$ is
\[\|u\|_{C^{k,\sigma }(V)}=\sup_{|\alpha |\le k,x\in
V}|D^\alpha u(x)|+\sup_{|\alpha |=k,x\ne y}|D^\alpha
u(x)-D^\alpha u(y)|\, |x-y|^{-\sigma }.)
\]
We then denote $k+\sigma =\tau $.

That a bounded domain $\Omega $ is $  C^{k,\sigma  }$
 means that there is an open cover 
$\{U_j\}_{j=1,\dots ,J}$ of
$\partial\Omega$ such that by an affine coordinate change for each $j$, $U_j$
is a box $\{\max_{k\le n}|y_k|<a_j\}$, and 
\[
\aligned
\Omega \cap U_j&=\{(y',y_n)\mid \max_{k<n}|y_k|<a_j,\,f_j(y')<y_n<a_j\},\\
\partial\Omega \cap U_j&=\{(y',y_n)\mid \max_{k<n}|y_k|<a_j,\,y_n=f_j(y')\},
\endaligned
\]
with $  C^{k,\sigma }$-functions $f _j$ such that $|f_j(y')|<a_j$ for $ \max_{k<n}|y_k|<a_j$. 
The diffeomorphism (coordinate change) 
\begin{equation}
F _j:(y',y_n)\mapsto
(y',y_n-f_j(y'))\label{tag4.1a}
\end{equation}
 is
then also $  C^{k,\sigma }$. The sets $U_j$
must be supplied with a suitable bounded  open set $U_0$ with closure contained
in $\Omega $, to get a full cover of $\comega$. 

For 
exterior domains, we cover $\partial\Omega $ similarly, then this must 
be supplied with a suitable open set $U_0$ with closure contained in $\Omega $
to get a full cover of $\comega$; here $U_0$
contains the complement of a ball, $U_0\supset {\Bbb R}^n\setminus
B(0,R')$.

For a perturbed halfspace, 
 we cover $\partial\Omega \cap B(0,R+1)$ as above, and supply this with
 $U_0=\{x\mid x_n>-\varepsilon , |x|>R\}$  to get a full cover of $\comega$.

The boundary $\partial\Omega $ will be denoted $\Sigma $.
We assume in the present section that $\Omega $ is $C^\infty $; then
$\Sigma $ is an $(n-1)$-dimensional $C^\infty $ manifold without
boundary. 

Let $A=\sum_{|\alpha |\le 2}a_\alpha D^\alpha$
with $C^\infty $ coefficients $a_\alpha $ given on a neighborhood
$\widetilde\Omega $ of $\comega$ (containing $U_0$ in the perturbed
halfspace case), and uniformly strongly
elliptic:
\[ \operatorname{Re}\sum_{|\alpha |= 2}a_\alpha(x) \xi ^\alpha\ge
c_0|\xi |^2,\text{ all }x\in \widetilde\Omega , \xi \in{\Bbb R}^n,
\]
$c_0>0$. The formal adjoint $A'=\sum_{|\alpha |\le 2}D^\alpha \bar
a_\alpha =\sum_{|\alpha |\le 2}a'_\alpha D^\alpha $ likewise has
$C^\infty $ coefficients  $a'_\alpha
$ and is strongly elliptic on $\widetilde \Omega $. We asume that the
coefficients and all their derivatives are bounded.

 We denote by $\Ama$ resp.\ $\Ami$ the maximal
resp.\ minimal realisations of $A$ in $L_2(\Omega )=H$; they act like $A$ in the
distribution sense and have the domains
\begin{equation}D(\Ama)=\{u\in L_2(\Omega )
\mid Au\in L_2(\Omega )\},\quad D(\Ami)=H^2_0(\Omega )\label{tag2.2}
\end{equation} 
(using $L_2$ Sobolev spaces). Similarly,  
$\Ama'$ and $\Ami'$ denote the maximal and minimal
realisations in $L_2(\Omega )$ of the formal adjoint $A'$; here
$\Ama={\Ami'}^*$, $\Ama'={\Ami}^*$.

Denote $\gamma _ju=(\partial_n ^ju)|_{\Sigma }$, where $\partial_n $ is the
derivative along the interior normal $\vec n $ at $\Sigma $. Let $s_0(x')$
be the coefficient of $-\partial_n ^2$ when $A$ is written in terms
of normal and tangential derivatives at $x'\in \Sigma $; it is bounded
with bounded inverse. Denoting 
\begin{equation}
s_0\gamma _1=\nu _1,\quad \bar s_0\gamma _1=\nu _1',\label{tag2.1a}
\end{equation}
we have the Green's formula for $A$
valid for $u,v\in H^2(\Omega )$, 
\begin{equation}
(Au,v)_{L_2(\Omega )}-(u,A'v)_{L_2(\Omega )}=(\nu  _1u,\gamma _0v)_{L_2(\Sigma )}-
(\gamma _0u,\nu '_1v+\mathcal A'_0\gamma _0v)_{L_2(\Sigma )},\label{tag2.1}
\end{equation}
where $\mathcal A'_0$ is  a certain first-order differential
operator over $\Sigma $. The formula extends e.g.\ to $u\in H^2(\Omega )$, 
$v\in D(\Ama')$, as
\begin{equation}
(Au,v)_{L_2(\Omega )}-(u,A'v)_{L_2(\Omega )}=(\nu  _1u,\gamma _0v)_{\frac12,-\frac12}-
(\gamma _0u,\nu '_1v+\mathcal A'_0\gamma _0v)_{\frac32,-\frac32},\label{tag2.1b}
\end{equation}
where $(\cdot,\cdot)_{s ,-s }$ denotes the duality pairing
between $H^s (\Sigma )$ and $H^{-s }(\Sigma
)$. (Cf.\ Lions and Magenes \cite{LM68} for this and the next results.)

The Dirichlet realisation $A_\gamma $ is defined as usual by
variational theory (the Lax-Milgram lemma); it is the restriction of $\Ama$ with
domain
\[
D(A_\gamma )=D(\Ama)\cap H^1_0(\Omega)=H^2(\Omega )\cap H^1_0(\Omega),
\]
where the last equality follows by elliptic regularity theory. By addition
of a constant to $A$ if necessary, we can assume that the spectrum of
$A_\gamma $ is contained in $\{\lambda \in {\Bbb C}\mid \operatorname{Re}\lambda >0\}$.
For $\lambda \in \varrho (A_\gamma )$, $s\in{\Bbb R}$, 
let 
\begin{equation}
Z^s_\lambda (A)=\{u\in H^s(\Omega )\mid (A-\lambda )u=0\};\label{tag2.3}
\end{equation}
it is a closed subspace of $H^s(\Omega )$. The trace operators $\gamma
_0$, $\gamma _1$ and $\nu  _1$ extend by
continuity to continuous maps
\begin{equation}
\gamma _0: Z^s_\lambda (A)\to H^{s-\frac12}(\Sigma ),\quad
\gamma _1,\nu _1: Z^s_\lambda (A)\to H^{s-\frac32}(\Sigma ),\label{tag2.4}
\end{equation}
for all $s\in{\Bbb R}$. When $\lambda \in \varrho (A_\gamma )$, let
$K^\lambda _{\gamma }:\varphi \mapsto u$ denote the Poisson operator
from $H^{s-\frac12}(\Sigma )$ to $H^s(\Omega )$ solving the
semi-homogeneous Dirichlet problem
\begin{equation}
(A-\lambda )u=0\text{ in }\Omega ,\quad \gamma _0u= \varphi \text{ on }\Sigma .\label{tag2.5}
\end{equation}
It is well-known that $K^\lambda _\gamma $ maps homeomorphically 
\begin{equation}
K^\lambda _\gamma :H^{s-\frac12}(\Sigma )\simto Z^s_\lambda (A),\label{tag2.6}
\end{equation}
for all $s\in {\Bbb R}$, with $\gamma _0$ acting as an inverse
there. The analogous operator for $A'-\bar\lambda $ is denoted
$K^{\prime\bar\lambda }_\gamma $. 

We shall now recall from \cite{BGW08, G68} how the statements in
Section \ref{Section1} are interpreted in terms of boundary
conditions. In the rest of this section, we abbreviate $H^s(\Sigma )$
to $H^s$.
With the notation from Section 1, 
\begin{equation}
Z^0_0(A)=Z,\quad Z^0_0(A')=Z',\quad Z^0_\lambda (A )=Z_\lambda ,\quad
Z^0_\lambda (A' )
=Z'_{\lambda }.\label{tag2.6a} 
\end{equation}
We denote by $\gamma _{Z_\lambda }$ the restriction of $\gamma _0$ to a mapping from
$Z_\lambda$ (closed subspace of $L_2(\Omega )$) to 
$H^{-\frac12} $; its adjoint $\gamma
_{Z_\lambda }^*$ goes from $H^{\frac12} $ to $Z_\lambda $:
\begin{equation}
\gamma _{Z_\lambda }: Z_\lambda \simto H^{-\frac12} ,\text{
with adjoint }\gamma _{Z_\lambda }^*:  H^{\frac12} \simto
Z_\lambda .\label{tag2.6c}
\end{equation}
There is a similar notation for the primed operators. When $\lambda
=0$, this index is left out.

These homeomorphisms allow ``translating'' an operator $T:Z\to Z'$ to
an operator $L:H^{-\frac12} \to H^\frac12 $, as in the diagram
\begin{equation}
\CD
Z     @>  \gamma _Z  >>    H^{-\frac12}\\
@VTVV           @VV  L  V\\
   Z'  @>>(\gamma _{Z'}^*)^{-1} >   H^{\frac12}\endCD \hskip1cm
   D(L)=\gamma _0D(T),
\label{tag2.6b}
 \end{equation}
whereby $(Tz,z')=(L\gamma _0z,\gamma _{0}z')_{\frac12,-\frac12}$.

We moreover define the Dirichlet-to-Neumann operators for each $\lambda \in
\varrho (A_\gamma )$, 
\begin{equation}
P^\lambda _{\gamma _0,\nu  _1}=\nu  _1K^\lambda _\gamma ;\quad 
P^{\prime\bar\lambda }_{\gamma _0,\nu ' _1}=\nu ' _1K^{\prime\bar\lambda }_\gamma ;
 \label{tag2.7}\end{equation}
they are  first-order elliptic pseudodifferential operators over
$\Sigma $, continuous from
$H^{s-\frac12} $ to $H^{s-\frac32} $
for all $s\in{\Bbb R}$, and Fredholm in case $\Sigma $ is
bounded. (Their pseudodifferential nature and ellipticity was explained e.g.\ in \cite{G71}).

For general trace maps $\beta $ and $\eta $ we write
\begin{equation}
P^\lambda _{\beta ,\eta }: \beta u\mapsto \eta u,\quad u\in
Z^s_\lambda (A),\label{tag2.8a}\end{equation}
when this operator is well-defined.

Introduce the trace operators $\Gamma $ and $\Gamma '$ (from \cite{G68},
where they were called $M$ and $M'$)  by 
\begin{equation}
\Gamma =\nu _1-P^0_{\gamma _0,\nu _1}\gamma _0=\nu _1A_\gamma ^{-1}\Ama,\quad \Gamma'
=\nu ' _1-P^{\prime 0}_{\gamma _0,\nu ' _1}\gamma _0=\nu _1'(A^*_\gamma )^{-1}\Ama'.\label{tag2.9}
\end{equation}
Here $\Gamma $ and $\Gamma '$ map $D(\Ama)$ resp.\ $D(\Ama')$ continuously onto
$H^{\frac12} $.
With these pseudodifferential boundary operators there is a generalized Green's formula valid {\it for all} 
$u\in D(\Ama)$, $v\in D(\Ama')$:
\begin{equation}
(Au,v)_{L_2(\Omega )}-(u,A'v)_{L_2(\Omega )}=(\Gamma  u,\gamma _0v)_{\frac12,
-\frac12}-(\gamma _0u, \Gamma 'v)_{-\frac12,\frac12}.\label{tag2.10}
\end{equation}
In particular,
\begin{equation}
(Au,w)=(\Gamma u,\gamma _0w)_{\frac12,-\frac12}\text{ for all
}w\in Z^0_0(A')=Z'.\label{tag2.11}
\end{equation}
(Cf.\ \cite{G68}, Th. III 1.2.) By composition with suitable isometries
$\Lambda _{t}:H^s(\Sigma )\to H^{s-t}(\Sigma )$, \eqref{tag2.10} can be
turned into a standard boundary triplet formula
\begin{equation}
(Au,v)_{L_2(\Omega )}-(u,A'v)_{L_2(\Omega )}=(\Gamma _1 u,\Gamma '_0v)_{L_2(\Sigma )}-(\Gamma _0u, \Gamma '_1v)_{L_2(\Sigma )},\label{tag2.10a}
\end{equation}
with $\Gamma _1=\Lambda _\frac12 \Gamma $, $\Gamma ' _1=\Lambda
_\frac12 \Gamma' $, $\Gamma _0=\Gamma '_0=\Lambda _{-\frac12}\gamma
_0$ and $\mathcal H=\mathcal K=L_2(\Sigma )$.

There is a general ``translation'' of the abstract results in Section
1 to statements on closed realisations $\wA$ of $A$. First let $\wA$
correspond to $T:Z\to Z'$ (i.e., assume $V=Z$, $W=Z'$). Then in view of
\eqref{tag2.6b} and \eqref{tag2.11}, the defining equation in Theorem
\ref{Theorem1.1} is turned into
\[
(\Gamma  u,\gamma _0z')_{\frac12,-\frac12}=(L\gamma _0u,\gamma _0z')_{\frac12,-\frac12},\text{ all }z'\in Z'.
\]
Since $\gamma _0z'$ runs through $H^{-\frac12}$, this means that
$\Gamma u=L\gamma _0u$,
also written 
\begin{equation}
\nu _1u=(L+P^0_{\gamma _0,\nu _1})\gamma _0u.\label{tag2.11a}
\end{equation}
Thus $\wA$ represents a
{\it Neumann-type condition} 
\begin{equation}
\nu  _1u=C\gamma _0u,\text{ with }C=L+P^0_{\gamma _0,\nu  _1}.\label{tag2.11b} 
\end{equation}

This allows  all
first-order $\psi $do's $C$ to enter, namely by letting $L$ act as
$C-P^0_{\gamma _0,\nu  _1}$. 

{\it The elliptic case:}
Consider a Neumann-type boundary condition
\begin{equation}
\nu _1u=C\gamma _0u,\label{tag2.12}
\end{equation}
where $C$ is a first-order classical $\psi $do on $\Sigma $. Let
$\wA$ be the restriction of $\Ama$ with domain
\begin{equation}
D(\wA)=\{u\in D(\Ama)\mid \nu _1u=C\gamma _0u\}.\label{tag2.13}
\end{equation}
Now the boundary condition satifies the
Shapiro-Lopatinski\u\i{} condition (is {\it elliptic}) if and only if
$L$ is elliptic; then in fact
\begin{equation}
D(\wA)=\{u\in H^2(\Omega )\mid \nu _1u=C\gamma _0u\}.\label{tag2.14}
\end{equation}
Then the adjoint $\wA^*$ equals
the operator that is defined
similarly from $A'$ by the boundary condition 
\begin{equation}
\nu' _1v=(C^*-\mathcal A'_0)\gamma _0v,\label{tag2.15}
\end{equation}
likewise elliptic.

When we do the above considerations  for $\wA-\lambda $, we get
$L^\lambda $ satisfying the diagram
%\vskip.2cm
\[
\CD
Z     @>{E^\lambda _{Z}} >>Z _\lambda     @>  \gamma _{Z_\lambda }  >>    H^{-\frac12}\\
@V  T+G^\lambda _{Z,Z'} VV @VT^\lambda VV           @VV  L^\lambda   V\\
 Z' @>>(F^{\prime\bar\lambda }_{Z'})^* >  Z'_{\bar\lambda }
@>>(\gamma _{Z'_{\bar\lambda }}^*)^{-1} >   H^{\frac12}\endCD
\hskip.5cm D(L^\lambda )=D(L).
 \] 
Here the horizontal maps are homeomorphisms, and they compose as $\gamma _{Z_\lambda }E^\lambda _Z=\gamma _Z$, $(\gamma _{Z'_{\bar\lambda }}^*)^{-1}(F^{\prime\bar\lambda }_{Z'})^*=(\gamma _{Z'}^*)^{-1}$,
so
\[   L^\lambda =\gamma _Z^{-1}(T+G^\lambda _{Z,Z'})\gamma _{Z'}^*. 
\]
In
terms of $L^\lambda $,
the boundary condition reads:
\[
\nu _1u=(L^\lambda +P^\lambda _{\gamma _0,\nu _1})\gamma _0u.
\]
Note that $L^\lambda +P^\lambda _{\gamma _0,\nu _1}=C=L +P^0
_{\gamma _0,\nu _1}$, so 
\[
 L^\lambda =L+P^0 _{\gamma _0,\nu _1}-P^\lambda _{\gamma _0,\nu _1}.
\]

As shown in \cite{BGW08}, this leads to:

\begin{theorem}\label{Theorem2.1} Assumptions as in the start of
  Section {\rm \ref{Section2}}, with $C^\infty $ domain and operator.
 Let $\wA$ correspond to $T:Z\to Z'$,
carried over to $L:H^{-\frac12}\to H^{\frac12}$. Then $\wA$ represents
the boundary condition {\rm \eqref{tag2.11b}}. Moreover:

{\rm (i)} For $\lambda \in \varrho (A_\gamma )$,  $P^0 _{\gamma
_0,\nu _1}-P^\lambda _{\gamma _0,\nu _1}\in \mathcal L(H^{-\frac12},
H^{\frac12})$ and
\[
L^\lambda =L+P^0 _{\gamma _0,\nu _1}-P^\lambda _{\gamma _0,\nu _1}.
\]

{\rm (ii)} For $\lambda \in \varrho (\wA)$, there is a related
$M$-function $\in\mathcal L(H^{\frac12}, H^{-\frac12})$
 \[
M_L(\lambda ) =\gamma _0\bigl(I-(\wA-\lambda )^{-1}(\Ama-\lambda
)\bigr)A_\gamma  ^{-1}\inj_{Z'\to H}\gamma _{Z'}^*.
\]

{\rm (iii)} For $\lambda \in \varrho (\wA)\cap \varrho (A_\gamma )$,
\[
M_L(\lambda )=-(L+P^0 _{\gamma _0,\nu _1}-P^\lambda _{\gamma
_0,\nu _1})^{-1}=-(L^\lambda )^{-1}.
\]

{\rm (iv)} For $\lambda \in \varrho (A_\gamma )$, 
\begin{align}
\operatorname{ker}(\wA-\lambda )&=K^\lambda _\gamma
\operatorname{ker}L^\lambda,\nonumber\\
\operatorname{ran}(\wA-\lambda )&=\gamma _{Z'_{\bar\lambda }}^*
\operatorname{ran}L^\lambda+\operatorname{ran}(\Ami-\lambda ),
 \label{tag2.17}
\end{align}
 so that $H\setminus (\operatorname{ran}(\wA-\lambda ))=Z^{\prime}_{\bar\lambda
}\setminus (\gamma _{Z'_{\bar\lambda }}^*
\operatorname{ran}L^\lambda)$.

{\rm (v)} For $\lambda \in \varrho (\wA)\cap \varrho (A_\gamma )$
there is a
Kre\u\i{}n resolvent formula:
\begin{align} (\wA-\lambda )^{-1}&=(A_\gamma -\lambda )^{-1}-\inj_{Z_\lambda \to H}\gamma _{Z_\lambda }^{-1} 
M_L(\lambda )(\gamma _{Z'_{\bar\lambda }}^*)^{-1}\pr_{Z'_{\bar\lambda }}\nonumber\\
&=(A_\gamma -\lambda )^{-1}-K^\lambda _{\gamma } 
M_L(\lambda )(K^{\prime\bar\lambda }_\gamma )^*
.\label{tag2.19}\end{align}

{\rm (vi)} In particular, if $C$ is a $\psi $do of order $1$ such that
$C-P^0_{\gamma _0,\nu _1}$ is elliptic, and $ \varrho (\wA)\cap
\varrho (A_\gamma )\ne \emptyset$, then $D(L)=H^{\frac32} $,
and 
\begin{equation}
M_L(\lambda )=-(C-P^\lambda _{\gamma
_0,\nu _1})^{-1}
\label{tag2.18a}
\end{equation} is elliptic of order $-1$ for all $\lambda \in \varrho (\wA)$.
Here $\wA$ satisfies {\rm \eqref{tag2.14}} with {\rm \eqref{tag2.11b}}.
\end{theorem}

Note that with the notation \eqref{tag2.8a}, $C-P^\lambda _{\gamma
_0,\nu _1}=-P^\lambda _{\gamma _0,\nu _1-C\gamma _0}$, and $M_L(\lambda )=P^\lambda _{\nu _1-C\gamma _0,\gamma _0}$.

Observe the simple last formula in \eqref{tag2.19}, where $K^\lambda
_\gamma $ is the Poisson operator for $A-\lambda $, the 
adjoint being a trace operator of class
zero.

The Kre\u\i{}n formula is consistent with formulas found for selfadjoint
cases with Robin-type conditions in other works, such as 
Posilicano  \cite{P08}, Posilicano and Raimondi  \cite{PR08}, Gesztesy and
Mitrea \cite{GM08}, when one observes that
\begin{equation}
(K^{\prime\bar\lambda }_\gamma )^*=\nu _1(A_\gamma -\lambda )^{-1};\label{tag2.20}
\end{equation}
 this follows from the fact that for $\varphi \in H^{-\frac12}(\Sigma
)$ and $v=K^{\prime\bar\lambda} _\gamma \varphi $, $f\in L_2(\Omega
)$ and $u=(A_\gamma -\lambda )^{-1}f$, one has using Green's formula \eqref{tag2.1b}: 
\begin{multline}
(f,K^{\prime\bar\lambda }_\gamma \varphi )_{L_2(\Omega )}=
((A-\lambda )u,v)_{L_2(\Omega )}-(u,(A'-\bar\lambda )v)_{L_2(\Omega )}\nonumber\\
=(\nu _1u,\gamma
_0v)_{\frac12,-\frac12}-(\gamma _0u,\nu '_1v+\mathcal A'_0\gamma _0v)_{\frac32,-\frac32}=(\nu _1(A_\gamma -\lambda )^{-1}f,\varphi )_{\frac12,-\frac12}.\nonumber
\end{multline}

For the general case of $\wA$ corresponding to $T:V\to W$ with
subspaces  $V\subset Z$, $W\subset Z'$, there is a related 
``translation'' to boundary
conditions. Details are given in \cite{BGW08}, let us here just
mention some ingredients:

We use the notation in \eqref{tag2.10a}\,ff.
 Set 
\[X_1=\overline{\Gamma _0 D(\wA)}=\Lambda _{-\frac12}\gamma _0V\subset L_2(\Sigma ) ,\quad Y_1=\overline{\Gamma
_0 D(\wA^*)}=\Lambda _{-\frac12}\gamma _0W\subset L_2(\Sigma ),\] 
where $\Gamma _0$ restricts to homeomorphisms
\[
\Gamma _{0,V}:V\simto X_1,\quad  \Gamma _{0,W}:W\simto Y_1.
\] 
Then $T:V\to W$ is carried over to $L_1:X_1\to Y_1$ by
\[
\CD
V     @>  \Gamma _{0,V}  >>    X_1\\
@VTVV           @VV  L_1  V\\
   W  @>>(\Gamma _{0,W}^*)^{-1} >   Y_1 \endCD \hskip1cm D(L_1)=\Gamma _0D(T),
 \]
The boundary condition is:
\[
\Gamma _0u\in D(L_1), \quad L_1\Gamma _0u=\pr_{Y_1}\Gamma _1u.
\]
There is a similar reduction for $\wA-\lambda $ when
$\lambda \in \varrho (A_\gamma )$, and we 
find that
\[
L_1^\lambda =L_1+\pr_{Y_1}\Lambda _\frac12(P^0_{\gamma _0,\nu
_1}-P^\lambda _{\gamma _0,\nu _1})\Lambda _\frac12 \inj_{X_1\to L_2(\Sigma )}.
\]
There is an $M$-function $M_{L_1}(\lambda ):Y_1\to X_1$ defined for
$\lambda \in \varrho (\wA)$. It equals $-(L^\lambda _1)^{-1}$ when
$\lambda \in \varrho (\wA)\cap \varrho (A_\gamma )$, and there is then
a Kre\u\i{}n resolvent formula
\begin{align}
(\wA-\lambda )^{-1}&=(A_\gamma -\lambda )^{-1}-\inj_{V_\lambda \to
  H}\Gamma _{0,V_\lambda }^{-1} M_{L_1}(\lambda )(\Gamma
_{0,W_{\bar\lambda }}^*)^{-1}\pr_{W_{\bar\lambda }} \nonumber\\
&=(A_\gamma -\lambda )^{-1}-K^\lambda _{\gamma ,X_1} M_{L_1}(\lambda )(K^{\prime\bar\lambda }_{\gamma ,Y_1})^*;
\nonumber\end{align}
here $K^\lambda _{\gamma ,X_1}: X_1\subset L_2(\Sigma )
\overset{\Lambda _\frac12} \longrightarrow H^{-\frac12}(\Sigma ) \overset
{K^\lambda _\gamma }\longrightarrow L_2(\Omega )$.

For higher order elliptic operators, and systems, there are similar
results on $M$-functions and Kre\u\i{}n resolvent formulas, see
\cite{BGW08}. In such cases there occur interesting subspace
situations where $X$ and $Y$ are (homeomorphic to) full products of
Sobolev spaces over $\Sigma $.

\medskip

\section{The nonsmooth $\psi $dbo calculus.}\label{Section3}

The study of the smooth case was formulated in \cite{BGW08} in terms of the
pseudodifferential boundary operator ($\psi $dbo) calculus, which was initiated
by Boutet de Monvel \cite{B71} and further developed e.g.\ in
Grubb \cite{G84}, \cite{G96} (we refer to these works or to \cite{G09}
for details on the calculus). The $\psi $dbo
theory has been adapted to nonsmooth situations by Abels in
\cite{A05}, by use of ideas from the
adaptation of $\psi $do's to nonsmooth cases by Kumano-go and Nagase
\cite{KN78}, Taylor \cite{T91}. The operators considered by Abels have
symbols that satisfy the usual estimates in the conormal variables
$\xi ',\xi ,\eta _n$, pointwise in the space variable $x$, but are only
of class $C^{k,\sigma }$ in $x$ (so that the symbol estimates hold with
respect to $ C^{k,\sigma }$-norm in $x$). 
(For $\tau =k+\sigma $ integer, one could replace $C^{k,\sigma }$ by the so-called Zygmund space 
$C^{\tau }=B^\tau _{\infty ,\infty }$, which is slightly larger, and gives the
scale of spaces slightly better interpolation properties, cf.\
Abels \cite{A03, A05b}, but we shall let that aspect lie.) We call $(k,\sigma ) $ the
H\"older smoothness of the operator and its symbol.

The theory allows the operators to act beween $L_p$-based Besov and 
Bessel-potential spaces ($1<p<\infty $), but we shall here just use
it in the case $p=2$ (although an extension to $p\ne 2$ would
also be interesting). Some important 
results of \cite{A05} are:

%\end{document}

\begin{theorem}\label{Theorem3.1}

$1^\circ$ One has that
\begin{equation}
\mathcal A=\begin{pmatrix} P_++G& K\\ \quad&\quad\\T& S\end{pmatrix} : \begin{matrix}
H^{s+m}(\rnp)^N\\ \times \\ H^{s+m-\frac12}({\Bbb
R}^{n-1})^M\end{matrix} \to
\begin{matrix}
H^{s}(\rnp)^{N'}\\ \times \\ H^{s-\frac12}({\Bbb R}^{n-1})^{M'}\end{matrix} 
\label{tag3.2}\end{equation}
holds when $\mathcal A$ is a Green operator on $\rnp$ of order $m\in{\Bbb Z}$ and class $r$,
with H\"older smoothness $(k,\sigma )$, provided that (with $\tau =k+\sigma $)
\begin{enumerate}
\item  $
|s|<\tau \text{ if }N'\ne 0$, 
\item  $ |s-\tfrac12|<\tau \text{ if }M'\ne 0$,
\item  $s+m>r-\tfrac12\text{ if }N\ne 0$ (class restriction).
\end{enumerate}

$2^\circ$ Let $\mathcal A_1$ and $\mathcal A_2$ be as in $1^\circ$, with
symbols $a_1$ resp.\ $a_2$ and constants $k_1,\sigma _1,\tau _1,m_1,N_1,\dots$ resp.\ $k_2,\sigma _2,\tau _2,m_2,N_2,\dots$. Assume that
$N_2'=N_1, M_2'=M_1$, so that the operators can be composed. Let
 $k _3=\min \{k _1,k _2\}$,  $\sigma _3=\min \{\sigma  _1,\sigma _2\}$,  $\tau _3=\min \{\tau _1,\tau _2\}$, $0<\theta <\min \{1,\tau _2\}$.
The boundary
symbol composition $a_1\circ_na_2$ is a Green symbol $a_3$ of order
$m_3=m_1+m_2$, class $r_3=\max\{r_1+m_2,r_2\}$ and H\"older smoothness
$(k_3, \sigma _3)$, defining a Green operator $\mathcal A_3$. The remainder is
continuous:
\begin{equation}
\mathcal A_1\mathcal A_2-\mathcal A_3 : \begin{matrix}
H^{s+m_3-\theta }(\rnp)^{N_2}\\ \times \\ H^{s+m_3-\frac12-\theta }({\Bbb
R}^{n-1})^{M_2}\end{matrix} \to
\begin{matrix}
H^{s}(\rnp)^{N'_1}\\ \times \\ H^{s-\frac12}({\Bbb R}^{n-1})^{M'_1}\end{matrix} ,
\label{tag3.3}
\end{equation}
if the following conditions are satisfied:
\begin{enumerate}
\item  $|s|<\tau _3$ and $s-\theta >-\tau _2$ if $N'_1> 0$,
$|s-\frac12|<\tau _3$ and $s-\frac12-\theta >-\tau _2$ if $M'_1> 0$;
\item  $-\tau _2+\theta <s+m_1<\tau _2$ if $N_1> 0$, $-\tau _2+\theta
<s+m_1-\frac12<\tau _2$ if $M_1> 0$;
\item  $s+{m_1} > r_1-\frac12$ if $N_1> 0$, $s+m_3-\theta >r_2-\frac12$ if
$N_2> 0$ (class restrictions).
\end{enumerate}

$3^\circ$ Let $\mathcal A$ be as in $1^\circ$, and polyhomogeneous and
uniformly elliptic with principal symbol $a^0$ (here $N=N'>0$). Then there is
a Green operator $\mathcal B^0$ (the operator with symbol $(a^0)^{-1}$
if $m=0$)
of order $-m$, class $r-m$ and H\"older smoothness $(k,\sigma )
$, continuous in the opposite direction of $\mathcal A$, such that  
$\mathcal R=\mathcal
A\mathcal B^0-I$ is continuous:
\begin{equation}
\mathcal R: \begin{matrix}
H^{s-\theta }(\rnp)^N\\ \times \\ H^{s-\theta -\frac12}({\Bbb
R}^{n-1})^{M'}\end{matrix} \to
\begin{matrix}
H^{s}(\rnp)^{N}\\ \times \\ H^{s-\frac12}({\Bbb
  R}^{n-1})^{M}\end{matrix} ,
\label{tag3.4}
\end{equation}
if, with $\tau =k+\sigma $,
\begin{enumerate}
\item $
-\tau +\theta <s<\tau $;
\item $ s-\tfrac12 >-\tau +\theta $ if $M$ or $M'>0$;
\item $ s-\theta > r-m-\tfrac12$ (class restriction).
\end{enumerate}
\end{theorem}

See \cite{A05} (Theorems 1.1, 1.2 and 6.4). 
For integer $\tau $, the
results are worked out there for symbols in Zygmund spaces, but they
imply the results with H\"o{lder spaces, see also \cite{A03, A05b}. The class restrictions are
imposed even when the operators have $C^\infty $ coefficients. $\mathcal B^0$ is called a
parametrix of $\mathcal A$.

Abels has also generalized the calculus of
\cite{G96} for symbols depending on a parameter $\mu $ to nonsmooth coefficients; again the estimates in the
cotangent variables $\xi ',\xi ,\eta _n, \mu $ are the usual ones, but
valid in $x$  w.r.t.\ H\"older norms. 

We recall from the theory of $\psi $do's that $P$ is said to be ``in
$x$-form'' resp.\ ``in $y$-form'', when it is defined from a symbol
$p$ by
\[ Pu=c\int e^{ix\cdot \xi }p(x,\xi )\hat u(\xi
)\,d\xi ,\text{ resp.\ } Pu=c\int e^{i(x-y)\cdot \xi }p(y,\xi ) u(y )\,dyd\xi ,
\] 
$c=(2\pi )^{-n}$; the concept extends to $\psi $dbo's.
In Theorem \ref{Theorem3.1}, all the operators
labeled with $\mathcal A$ are in
$x$-form. 
%(cf.\ e.g.\ \cite{G96}, (2.1.33)\,ff.). 
So is $\mathcal
B^0$ when $m=0$; otherwise it is a composition of an operator in
$x$-form with an order-reducing operator system to the left, 
see Remark \ref{Remark3.2} below.
The adjoints of operators in $x$-form are operators in $y$-form. \cite{A05}
does not discuss the reduction from $y$-form to $x$-form; some
indications may be inferred from Taylor \cite{T00}, Ch.\ 1 §9. For
operators in $y$-form
one has at least the results that can be derived from the above
results by transposition.

\begin{remark}\label{Remark3.2} An important tool in the calculus is
``order-reducing
operators''. There are two types, one acting over the domain and one acting
over the boundary: 
\begin{align}
&\Lambda _{-,+}^r=\operatorname{OP}(\lambda ^r_-(\xi ))_+: H^t(\rnp)\simto H^{t-r}(\rnp),
\label{tag3.5}\\
&\Lambda _0^r=\operatorname{OP}'(\ang{\xi '}^r):H^t({\Bbb R}^{n-1})\simto H^{t-r}({\Bbb
R}^{n-1}),\text{ all $t\in {\Bbb
R}$,}\nonumber
\end{align}
 with inverses 
$\Lambda
_{-,+}^{-r}$ resp.\ $\Lambda _0^{-r}$. Here  $\lambda _-^r$ is the
``minus-symbol'' defined in
\cite{G90} Prop. 4.2 as a refinement of $(\ang{\xi '}-i\xi _n)^{r}$. 
In Theorem \ref{Theorem3.1} $3^\circ$, whereas $\mathcal B^0$ is the
operator with symbol $(a^0)^{-1}$ when $m=0$, one applies
the zero-order construction to $\mathcal A_1=\mathcal A\begin{pmatrix}\Lambda
  _{-,+}^{-m}&0\\0&\Lambda _0^{-m}\end{pmatrix}$ to define
$\mathcal B^0=\begin{pmatrix}\Lambda _{-,+}^{-m}&0\\0&\Lambda
  _0^{-m}\end{pmatrix}\mathcal B_1^0$ when $m\ne 0$.

It should be noted that when e.g.\ $P_+$ is as in Theorem
\ref{Theorem3.1} $1^\circ$, then 
\begin{equation}\Lambda _{-,+}^r P_+:H^{s+m}(\rnp)\to H^{s-r}(\rnp)
\text{ for }-\tau <s<\tau ,\label{tag3.6}\end{equation}
 whereas the
composition rule  Theorem \ref{Theorem3.1} $2^\circ$ shows that
$\Lambda _{-,+}^rP_+$ can be written as the sum of an operator in the
calculus $\operatorname{OP}'(\lambda ^r_{-,+}\circ_n p(x,\xi )_+)$ in $x$-form and a
remainder, such that the sum maps  $H^{s'+m+r}(\rnp)\to H^{s'}(\rnp)$ for $-\tau
<s'<\tau $; this gives a mapping property like in \eqref{tag3.6}
but with $-\tau +r<s<\tau +r$. This apparently extends the
range, but the decompositions into a primary part and a
remainder are not the same;
%It may seem convenient to have such
%different scales of spaces in the picture, but it should be noted 
%that the part of the operator
%lying in the calculus is not the same in the two representations; 
$\Lambda
_{-,+}^rP_+$ is not in $x$-form but is an operator in $x$-form composed to the left
with $\Lambda _{-,+}^r$, not equal to
$\operatorname{OP}'(\lambda ^r_{-,+}\circ_n p(x,\xi )_+)$. 
Compositions to the right with $\Lambda _{-,+}^r$ are simpler and
preserve $x$-form directly. We shall say that operators formed by
composing an operator in $x$-form with an order-reducing operator to
the left are ``in order-reduced $x$-form''.
\end{remark}

Coordinate changes give some inconveniences in the nonsmooth calculus because, in a $C^{k,\sigma }$-setting, the action of $D_j$ after a
$C^{k,\sigma }$-coordinate change gets Jacobian factors that are ${
C}^{k-1,\sigma }$, and higher powers $D^\alpha $ get coefficients in
$C^{k-|\alpha |,\sigma }$ (when $k-|\alpha |\ge 0$). 

We say that an operator is a generalized Green operator (of one of the
respective types) if it is the sum of an operator defined from 
symbols in the calculus and a remainder of lower order (for $s$ in an
interval, specified in each case or understood from the context).

\section{Resolvent formulas in the case of non-smooth domains.}\label{Section4}

To treat one difficulty at a time, we consider in the following the
case where the domain is nonsmooth, but the operator $A$ is given with
smooth coefficients (this
includes of course constant coefficients).

Let $\Omega $ be a open set in ${\Bbb R}^n$ of one of the three types
described in Section \ref{Section2}, of class $C^{k,\sigma }$. 
We still take $A$ with $C^\infty $-coefficients on a neighborhood
$\widetilde\Omega $ of $\comega$, as described in Section 2.

Recall from Grisvard \cite{G85} (Th. 1.3.3.1, 1.5.1.2, 1.4.1.1, 1.5.3.4)

\begin{theorem}\label{Theorem4.1} Let $\Omega $ be bounded and
  $ C^{k,\sigma }$, let $\tau =k+\sigma $.

$1^\circ$ When $\Phi $ is a $C^{k,\sigma }$-diffeomorphism, $\tau $
integer, then $u\in H^s_{\operatorname{loc}}\implies u\circ \Phi \in
H^s_{\operatorname{loc}}$ for $|s|\le \tau $.

$2^\circ$ One can for $|s|\le \tau $, integer,  define
$H^s(\Sigma  )$ to be the space of distributions $u$ on
$\Sigma  $ such that for each $j$, $ u\circ F_j^{-1}$ is in
$H^s$ on $\{y'\mid \max |y_k|\le a_j\}$. The trace
map $\gamma _0:H^s(\Omega )\to H^{s-\frac12}(\Sigma  )$ is
well-defined for $\frac12<s\le \tau $, and the trace
map $\gamma _1:H^s(\Omega )\to H^{s-\frac32}(\Sigma  )$ is
well-defined for $\frac32<s\le \tau $. There is a continous right
inverse of each map, and of the two maps jointly for $\frac32<s\le
\tau $.

$3^\circ$ Let $\varphi $ be $  C^{k_1,\sigma_1}$, $\tau
_1=k_1+\sigma _1$, then $u\mapsto \varphi u$ is
continuous in $H^s({\Bbb R}^n )$ for $| s|\le \tau _1$ if $\tau _1$ is
integer, $|s|<\tau _1$ if $\tau _1$ is non-integer. 

$4^\circ$ When $\tau \ge 2$ and $A$ is a
second-order differential operator on $\Omega $ in a divergence form
($A=-\sum_{j,k}\partial_ja_{jk }\partial_k+\sum_{k}a_k\partial_k+a_0$) with
$C^{0,1}$-coefficients, and we define the associated oblique Neumann
trace operators by 
\begin{equation}\nu _A=\sum_{j,k}n_ja_{jk} \gamma
_0\partial_k,\quad 
\nu _{A'}=\sum_{j,k}n_k\bar a_{jk} \gamma
_0\partial_j,\label{tag2.1e}
\end{equation}
there holds a  Green's formula 
\begin{multline}
(Au,v)_{L_2(\Omega )}-(u,A'v)_{L_2(\Omega )}\\=(\nu  _Au,\gamma _0v)_{L_2(\Sigma )}-
(\gamma _0u,\nu _{A'}v-{\sum}_kn_k\bar a_k\gamma _0v)_{L_2(\Sigma
  )},\label{tag2.1d} 
\end{multline}
for
$u,v\in H^2(\Omega )$. 

\end{theorem}

The Green's
formula \eqref{tag2.1d} can be
reorganized as \eqref{tag2.1}; for our $A$ with smooth coefficients,
$\nu _1$, $\nu _1'$ and $\mathcal A'_0$ get $C^{k-1,\sigma
}$-coefficients when $\Omega $ is $C^{k,\sigma }$.
 
We define the Dirichlet realisation $A_\gamma
$ of $A$, with domain $D(A_\gamma )=D(\Ama)\cap H^1_0(\Omega )$ by the
usual variational construction, and we shall assume that $A_\gamma $ is
invertible. Its adjoint is the analogous operator for $A'$. 

By the difference quotient method of Nirenberg \cite{N55} 
one has that $D(A_\gamma
)=H^2(\Omega )\cap H^1_0(\Omega )$ when $\tau \ge 2$ (this fact is
also derived below); detailed proofs are e.g.\ found in the
textbooks of Evans \cite{E98} (for $C^2$-domains) or McLean \cite{M00}
(for $C^{1,1}$-domains).

Also the extended Green's formula \eqref{tag2.1b} is valid when $\tau
\ge 2$;
this follows by an extension of the proof in
Lions and Magenes \cite{LM68}, as mentioned in \cite{G85} Remark
1.5.3.5. 
It follows that the generalized Green's formula \eqref{tag2.10} holds, when $\Gamma $ and
$\Gamma '$ are defined
by
\begin{equation}
\Gamma =\nu _1A_\gamma ^{-1}\Ama,\quad \Gamma'
=\nu _1'(A^*_\gamma )^{-1}\Ama'.\label{tag2.9a}
\end{equation}

\medskip

The local coordinates (cf.\ \eqref{tag4.1a}) are used to reduce the 
curved situation to the
flat situation; then the boundary becomes straight but nonsmoothness
is imposed on the symbols. 

In the following we work out what the nonsmooth $\psi $dbo method 
can give for the Dirichlet
problem; this can be regarded as a basic exercise in the
calculus (some other cases appear in works of Abels
and coauthors).

First we consider the case of a uniformly strongly elliptic
second-order operator on $\rnp$ --- which we for simplicity of notation also
call $A$ --- with H\"older smoothness $(k_1,\sigma _1)$ and $\tau _1 =k_1+\sigma _1$, together
with a Dirichlet trace operator,
\begin{equation}
\mathcal A=\begin{pmatrix} A\\\gamma _0\end{pmatrix} : 
H^{s+2}(\rnp) \to
\begin{matrix}
H^{s}(\rnp)\\ \times \\ H^{s+\frac32}({\Bbb R}^{n-1})\end{matrix};\label{tag4.2}
\end{equation}
it is continuous for
\begin{equation}
-\tau _1 <s<\tau _1 ,\quad s>-\tfrac32,\label{tag4.3}
\end{equation}
extended to $|s|\le \tau _1$ if integer (cf.\ Theorem \ref{Theorem4.1}
3$^\circ$). 
To prepare for an application of
 Theorem \ref{Theorem3.1}, we apply order-reducing operators (cf.\
 Remark \ref{Remark3.2}) to reduce to order 0, introducing
\begin{equation}
\mathcal A_1=\begin{pmatrix}I&0\\0&\Lambda _0^2\end{pmatrix}\mathcal A
\Lambda _{-,+}^{-2}=
\begin{pmatrix} A\Lambda _{-,+}^{-2}\\ \quad\\\Lambda _0^2\gamma _0\Lambda _{-,+}^{-2}\end{pmatrix} : 
H^{s}(\rnp) \to
\begin{matrix}
H^{s}(\rnp)\\ \times \\ H^{s-\frac12}({\Bbb R}^{n-1})\end{matrix}, \label{tag4.4}
\end{equation}
for $s$ as in \eqref{tag4.3}\,ff. By Theorem \ref{Theorem3.1}
$3^\circ$  it has a parametrix 
$\mathcal B^0_1$ of order 0 and class $-1$ defined from the
principal symbols, 
\begin{equation}
\mathcal B^0_1=\begin{pmatrix} R^0_1& K^0_1\end{pmatrix} :\begin{matrix}
H^{s}(\rnp)\\ \times \\ H^{s-\frac12}({\Bbb R}^{n-1})\end{matrix} \to
H^{s}(\rnp),\label{tag 4.5}\end{equation} 
for $s$ satisfying 
\begin{equation}
-\tau _1 +\tfrac12<s<\tau _1 ,\quad s>-\tfrac32;\label{tag4.5a}
\end{equation}
here the remainder $\mathcal R_1=\mathcal A_1\mathcal B^0_1-I$
satisfies
\begin{equation}
\mathcal R_1: \begin{matrix}
H^{s-\theta }(\rnp)\\ \times \\ H^{s-\theta -\frac12}({\Bbb
R}^{n-1})\end{matrix} \to
\begin{matrix}
H^{s}(\rnp)\\ \times \\ H^{s-\frac12}({\Bbb R}^{n-1})\end{matrix}
,\label{tag4.6} 
\end{equation}
when $0<\theta <\min\{1,\tau _1 \}$,
\begin{equation}
-\tau _1 +\tfrac12+\theta <s<\tau _1,\quad s > -\tfrac32+\theta .\label{tag4.7}
\end{equation}

Then the equation $\mathcal A_1\mathcal B_1^0=I+\mathcal R_1$, also written 
\[
\begin{pmatrix} I&0\\0&\Lambda _0^2\end{pmatrix}\mathcal A\Lambda _{-,+}^{-2}\mathcal B_1^0=I+\mathcal R_1,
\]
implies by composition to the left with $\begin{pmatrix} I&0\\0&\Lambda
  _0^{-2}\end{pmatrix}$ and to the right with $\begin{pmatrix} I&0\\0&\Lambda
    _0^2\end{pmatrix}$:
\begin{equation}
\mathcal A\Lambda _{-,+}^{-2}\mathcal B_1^0\begin{pmatrix} I&0\\0&\Lambda
  _0^2\end{pmatrix}=I+\mathcal R,\text{ with } \mathcal R=\begin{pmatrix} I&0\\0&\Lambda
  _0^{-2}\end{pmatrix}\mathcal R_1\begin{pmatrix} I&0\\0&\Lambda
    _0^2\end{pmatrix}.\nonumber
\end{equation}
Hence
\[ \mathcal B^0
=\Lambda _{-,+}^{-2}\mathcal B_1^0\begin{pmatrix} I&0\\0&\Lambda
  _0^2\end{pmatrix}=\begin{pmatrix} R^0& K^0\end{pmatrix}
\]
is a parametrix of $\mathcal A$, with
\begin{gather}
\mathcal A \mathcal B^0=I+\mathcal R,\label{tag4.12a}\\
\mathcal B^0 :\begin{matrix}
H^{s}(\rnp)\\ \times \\ H^{s+\frac32}({\Bbb R}^{n-1})\end{matrix} \to
H^{s+2}(\rnp),\quad\mathcal R: \begin{matrix}
H^{s-\theta }(\rnp)\\ \times \\ H^{s-\theta +\frac32}({\Bbb
R}^{n-1})\end{matrix} \to
\begin{matrix}
H^{s}(\rnp)\\ \times \\ H^{s+\frac32}({\Bbb R}^{n-1})\end{matrix}
,\label{tag4.6a} 
\end{gather}
for $s$ as in \eqref{tag4.5a} resp.\ \eqref{tag4.7}. With the notation
from Remark \ref{Remark3.2}, $\mathcal B^0$ is in order-reduced $x$-form.

\medskip

Now consider the situation where $A$ has smooth coefficients and the
domain is nonsmooth. 
We shall go through the parametrix and inverse
construction in the case where the H\"older smoothness of the domain
is $(1,1)$ so that $\tau =2$. We have the
direct operator
\begin{equation}
\mathcal A=\begin{pmatrix} A\\\gamma _0\end{pmatrix} : 
H^{s+2}(\Omega ) \to
\begin{matrix}
H^{s}(\Omega )\\ \times \\ H^{s+\frac32}(\Sigma )\end{matrix},\label{tag4.2}
\end{equation}
it is continuous for $-\frac32<s\le 0$
(recall the restriction $s+2\le 2$ coming from Theorem
\ref{Theorem4.1} $2^\circ$).

For each $i=1,\dots, J$, the diffeomorhism \eqref{tag4.1a} carries
$\Omega \cap U_j $ over to $V_j=\{(y',y_n)\mid
\max_{k<n}|y_k|<a_j,\,0<y_n<a_j-f_j(y')\}$, such that $\partial\Omega \cap U_j$ is
mapped to $\{(y',y_n)\mid \max_{k<n}|y_k|<a_j,\,y_n=0\}$.  When the smooth differential operator $A$ is transformed to local
coordinates in this way, the principal part of the
resulting operator $\underline A$ has H\"older smoothness 
$ (0,1)$, so here $\tau
_1=1$. 
In each of these charts 
one constructs a parametrix $\underline{\mathcal
  B}^0$ for $\begin{pmatrix}
\underline A\\\gamma _0\end{pmatrix}$ as above (the 
coefficients of $\underline A$ can
be assumed to be extended to $\crnp$). When $\Omega $ is bounded or is
an exterior domain, one uses for the set $U_0$
 a parametrix of $A$ without changing coordinates. In the perturbed
 halfspace case, for the set $U_0$ one extends $A$ smoothly to $\crnp$
and uses
a smooth version of  the above construction. These parametrices are carried
back to the curved situation and pieced together using a partition of
unity subordinate to the cover $\{U_0,U_1,\dots, U_J\}$, as indicated
in \cite{G96}, p.\ 228 (the first factor $\varphi _i$ in each term in (2.4.77)
should be replaced by a function $\eta _i\in C_0^\infty (U_i)$ such
that $\eta _i\varphi _i=\varphi _i$, to get preservation of the
principal symbol after summation). Here the coordinate changes
allow the smoothness to remain at $(0,1)$; cf.\ \cite{A05b}, in
particular Section 5.3 there. 
The sum over $i$ is then a
parametrix of \eqref{tag4.2}; its composition with $\mathcal A$ gives
the identity plus a remainder of lower order, for
values $s$ as indicated above.

In the subsequent compositions below, it will always be understood that
they take place in local coordinates (after decomposing the operators
in pieces supported in the $U_i$ by use of suitable partitions of
unity) and are taken back to the curved situation afterwards.

In the present construction, we shall actually
carry a spectral parameter along that
will be useful for discussions of invertibility. So we now
replace the originally given $A$ by $A-\lambda $, to be studied for
large negative $\lambda $.

The parametrix
will be of the form 
\begin{equation}
\mathcal B^0(\lambda )=\begin{pmatrix} R^0(\lambda )& K^0(\lambda )\end{pmatrix} :\begin{matrix}
H^{s}(\Omega )\\ \times \\ H^{s+\frac32}(\Sigma )\end{matrix} \to
H^{s+2}(\Omega );\label{tag4.14}\end{equation} 
with $(k_1,\sigma _1)=(0,1)$ the condition \eqref{tag4.5a} means that $-\frac12<s<1$, so that,
along with the restriction coming from Theorem \ref{Theorem4.1}, we
have altogether that
\begin{equation}
-\tfrac12<s\le 0 \label{tag4.14a}
\end{equation}
is allowed. The remainder maps as follows:
 \begin{equation}
\mathcal R(\lambda )=\mathcal A(\lambda )\mathcal B^0(\lambda )-I: \begin{matrix}
H^{s-\theta }(\Omega )\\ \times \\ H^{s-\theta +\frac32}(\Sigma )
\end{matrix} \to
\begin{matrix}
H^{s}(\Omega )\\ \times \\ H^{s+\frac32}(\Sigma )\end{matrix} \label{tag4.15}
\end{equation}
for 
\begin{equation}
-\tfrac12+\theta <s\le 0.\label{tag4.15a}
\end{equation}

In order to get hold of the exact inverse, we shall use an old trick
of Agmon \cite{A62}, which implies a useful
$\lambda $-dependent estimate of the remainder: Write $-\lambda =\mu
^2$ ($\mu >0$), introduce an extra variable $t\in S^1$, and replace
$\mu $ by $D_t=-i\partial_t$; let
\begin{equation}
\widehat A=A+D_t^2 \text{ on }\Omega \times S^1. \label{4.15b}
\end{equation}
Then $\widehat A$ is strongly elliptic on $\Omega \times S^1$, and by
the
preceding construction (carried out with local coordinates respecting
the product structure), 
\[
\widehat {\mathcal A} =\begin{pmatrix} \widehat A\\ \gamma
_0\end{pmatrix}\text{ has a parametrix }\widehat{\mathcal B}^0 ,
\]
with mapping properties of $\widehat{\mathcal B}^0 $ and the remainder
$\widehat{\mathcal R}=\widehat{\mathcal A} \widehat{\mathcal B}^0 -I$
as in 
\eqref{tag4.14} and \eqref{tag4.15}
with $\Omega ,\Sigma $ replaced by $\widehat\Omega =\Omega \times
S^1$, $\widehat\Sigma =\Sigma \times
S^1$.

For functions $w$ of the form $w(x,t)=u(x)e^{i\mu t}$, 
\[
\widehat {\mathcal A} w=\begin{pmatrix} ( A+\mu ^2)w\\ \gamma
_0w\end{pmatrix},
\]
and similarly, the parametrix $\widehat{\mathcal B}^0 $ and the remainder
$\widehat{\mathcal R}$ act on such functions like $\mathcal B^0
(\lambda )$ and
$\mathcal R(\lambda )$ applied in the $x$-coordinate.

Moreover, for  $w(x,t)=u(x)e^{i\mu t}$, $u\in \mathcal S({\Bbb R}^n)$,
\[
\|w\|_{H^s({\Bbb R}^n\times S^1)}\simeq \|(1-\Delta +\mu ^2
)^su(x)\|_{L_2({\Bbb R}^n)}\simeq \|(1+|\xi |^2+\mu ^2)^{s/2}\hat u(\xi )\|_{L_2},
\]
with similar relations for Sobolev spaces over other sets. Norms as in
the right-hand side are called $H^{s,\mu }$-norms; they were
extensively used \cite{G96}, see the Appendix there for the definition on subsets.
The important observation is now that when $s'<s$ and
$w(x,t)=u(x)e^{i\mu t}$, then
\begin{multline}
\|w\|_{H^{s'}({\Bbb R}^n\times S^1)}\simeq \|(1+|\xi |^2+\mu
^2)^{s'/2}\hat u(\xi )\|_{L_2}\nonumber \\
\le  \ang\mu ^{s'-s}\|(1+|\xi |^2+\mu ^2)^{s/2}\hat u(\xi )\|_{L_2}\simeq\ang\mu ^{s-s'}  \|w\|_{H^s({\Bbb R}^n\times S^1)},\nonumber
\end{multline}
with constants independent of $u$ and $\mu $. Analogous estimates hold with
${\Bbb R}^n$ replaced by $\Omega $ or $\Sigma $.

Applying this principle to the estimates of the remainder $\widehat{\mathcal R}$, we find 
that
\begin{multline}
\|\mathcal R(\lambda )u\|_{H^{s,\mu }(\Omega )\times H^{s+\frac32,\mu
}(\Sigma )}
\le c_s\|u\|_{H^{s-\theta ,\mu }(\Omega )\times
H^{s-\theta +\frac32,\mu }(\Sigma )}\nonumber
\\ \le c'_s\ang\mu ^{-\theta }\|u\|_{H^{s,\mu }(\Omega )\times
H^{s+\frac32,\mu }(\Sigma )}
\label{tag4.18}
\end{multline}
for $s$ as in \eqref{tag4.15a}.

For each $s$, take a fixed $\lambda $ with $|\lambda |$ so large that
$c'_s\ang\mu ^{-\theta }\le \frac12$. Then $I+\mathcal R(\lambda )$ has the
inverse $I+\mathcal R '(\lambda )=I+\sum_{k\ge 1}(-\mathcal R(\lambda ))^k$
(converging in the operator norm for operators on $H^{s,\mu }(\Omega
)\times H^{s+\frac32,\mu }(\Sigma )$), and 
\[
\mathcal A (\lambda )\mathcal B ^0(\lambda )(I+\mathcal R' (\lambda ))=I.
\]
This gives a right inverse
\[
\mathcal B (\lambda )=\mathcal B ^0(\lambda )+\mathcal B^0(\lambda )\mathcal R'(\lambda )=\begin{pmatrix}
R (\lambda )& K (\lambda )\end{pmatrix},
\]
with the same Sobolev space continuity as $\mathcal B ^0(\lambda )$, and
$\mathcal B ^0(\lambda )\mathcal R'(\lambda )$ of lower order.
 Since
\begin{equation}
\mathcal A (\lambda )\mathcal B (\lambda )=\begin{pmatrix} (A-\lambda )R (\lambda
)& (A-\lambda )K (\lambda )\\ \gamma _0R (\lambda
)& \gamma _0 K (\lambda )
\end{pmatrix}=\begin{pmatrix} I&0\\0&I\end{pmatrix},\label{tag4.19}
\end{equation}
$R (\lambda )$ solves
\begin{equation}
(A-\lambda )u=f,\quad \gamma _0u=0, \label{tag4.20}
\end{equation}
and $K (\lambda )$ solves
\begin{equation}
(A-\lambda )u=0,\quad \gamma _0u=\psi .\label{tag4.21}
\end{equation}

For such large $\lambda $, $R(\lambda )$ coincides with the resolvent
of $A_\gamma $ defined by variational theory, and $K(\lambda )$ is the
Poisson-type operator we called $K^\lambda _\gamma $ in Section \ref{Section2};
\begin{equation}
(A_\gamma -\lambda )^{-1}:H^s(\Omega )\to H^{s+2}(\Omega
),\quad K^\lambda _\gamma :H^{s+\frac32}(\Sigma )\to H^{s+2}(\Omega ),
\label{tag4.26}\end{equation}
for $s$ satisfying \eqref{tag4.14a}.

The mapping properties extend to all the $\lambda $ for which the
operators are well-defined, especially to $\lambda =0$.
For $A_\gamma ^{-1}$, this goes as follows: When $u\in H^1(\Omega )$ and $f\in
H^s(\Omega )$ with $s<1$, $f+\lambda u$ is likewise in $H^s(\Omega
)$. Then $A_\gamma u=f+\lambda u$ allows the conclusion $u\in
H^{s+2}(\Omega )$.  The argument works for all $s$ satisfying \eqref{tag4.14a} 
(for each such $s$, there is room to take $\theta >0$ so small that 
\eqref{tag4.15a} is satified. Moreover, since 
$A_\gamma ^{-1}-(A_\gamma -\lambda)^{-1
}=-\lambda A_\gamma ^{-1}(A_\gamma -\lambda )^{-1}$ is of lower order
than $A_\gamma ^{-1}$, $A_\gamma ^{-1}$ equals a nonsmooth $\psi $dbo
 plus a lower-order remainder. 

The
Poisson operator solving \eqref{tag4.21} can be further described as
follows (for all $\lambda \in \varrho (A_\gamma )$): There is a right inverse $\mathcal K:H^{s+\frac32}(\Sigma )\to
H^{s+2}(\Omega )$ of $\gamma _0$ for $-\frac32<s\le
0$ (cf.\ Theorem \ref{Theorem4.1} $2^\circ$). When we set $v=u-\mathcal K\varphi $,
we find that $v$ should solve
\[
(A-\lambda )v=-(A-\lambda )\mathcal K\varphi ,\quad \gamma _0v=0,
\]
to which we apply the preceding results; then when $\lambda \in \varrho (A_\gamma )$,
\begin{equation}
K^\lambda _\gamma =\mathcal K-(A_\gamma -\lambda )^{-1}(A-\lambda )\mathcal K;\label{tag4.21b}
\end{equation}
solves \eqref{tag4.21} uniquely.
It maps
$H^{s+\frac32}(\Sigma )\to H^{s+2}(\Omega )$ for $s$
satisfying \eqref{tag4.14a}.

Since our original operator had $C^\infty $
coefficients, the same construction works for the adjoint Dirichlet
problem, so we also here get the mapping properties
\begin{equation}
(A_\gamma ' -\bar\lambda )^{-1}:H^s(\Omega )\to H^{s+2}(\Omega ), \quad K^{\prime\bar\lambda }_\gamma
:H^{s+\frac32}(\Sigma )\to H^{s+2}(\Omega ),\label{tag4.21a}
\end{equation}
for $s$ satisfying \eqref{tag4.14a}.

The condition $s>-\frac12$ prevents the Poisson operator from starting
from $H^{-\frac12}(\Sigma )$, which would be needed for an 
analysis as in Section \ref{Section2}. Fortunately, it is possible to get
supplementing information in 
other ways.

By
\eqref{tag2.1b} we
have, analogously to \eqref{tag2.20}, that $K^\lambda _\gamma $ is the adjoint
of a trace operator of class 0 as follows:
\begin{equation}
K^{\lambda }_\gamma =(\nu '_1(A'_\gamma -\bar\lambda )^{-1})^*;\label{tag4.27}
\end{equation}
(it is used here that $\mathcal A'_0\gamma _0(A'_\gamma -\bar\lambda
)^{-1}=0$).

Now use the mapping property in \eqref{tag4.21a}.
The resolvent can be composed with $\nu '_1$ for $s>-\frac12$, so
\[
\nu _1'(A'_\gamma -\lambda )^{-1}=(K^\lambda _\gamma )^*:H^s(\Omega )\to H^{s+\frac12}(\Sigma )\text{
for }-\tfrac12<s\le 0.
\]
It follows that 
\begin{equation}
K^\lambda _\gamma : H^{s'-\frac12}(\Sigma )\to H^{s'}(\Omega
),
\label{tag4.29}
\end{equation}
when $0\le s'<\tfrac12$. In particular, $s'=0$ is allowed.

Taking this together with the larger values that were covered by \eqref{tag4.26},
we find that \eqref{tag4.29} holds for
\begin{equation}
0\le s'\le 2;\label{tag4.31}
\end{equation}
the intermediate values are included by interpolation. We denote $s'$
by $s$ from here on. 

One can analyze the structure of $K^\lambda _\gamma $ for the low
values of $s$ further, decomposing it into terms belonging to the
calculus and lower-order remainders. There is a difficulty here in the
fact that order-reducing operators as well as operators in $y$-form
enter, and both types affect the $s$-values for which the
decompositions and mapping
properties are valid (cf.\ Remark \ref{Remark3.2}). We refrain from including a deeper analysis.

There is a similar result for $K^{\prime\bar\lambda }_\gamma $. The
adjoints also extend, e.g.
\begin{equation}
(K^{\prime\bar\lambda }_\gamma)^*: H^s_0(\comega )\to H^{s+\frac12}(\Sigma
),\text{ for }-2\le s\le 0;\label{tag4.32}
\end{equation}
recall that $H^s_0(\comega)=H^s(\Omega )$ when $|s|< \frac12$. To sum
up, we have shown:

\begin{theorem}\label{Theorem4.2} When $\Omega $ is $ C^{1,1}$
  and $A$ has $C^\infty $-coefficients, the solution operators
$K^\lambda _\gamma $ and $K^{\prime\bar\lambda }_\gamma $ for {\rm
  \eqref{tag2.5}} and its primed version 
map $H^{s-\frac12}(\Sigma )$ to $H^s(\Omega
)$ for $0\le s\le 2$. They 
are generalized Poisson
operators in the sense that for $s\in \,]\frac32,2]$, they can be written as
the sum of a Poisson operator  of H\"older smoothness $(0,1)$, in
order-reduced $x$-form, and a
lower order operator.
\end{theorem}

The next step is to study $P^\lambda _{\gamma _0,\nu  _1}
=\nu  _1K^\lambda _\gamma $ and $ 
P^{\prime\bar\lambda }_{\gamma _0,\nu ' _1}=\nu '
_1K^{\prime\bar\lambda }_\gamma $, cf.\ \eqref{tag2.7}\,ff.

We have immediately from the mapping properties established above,
that
\begin{equation}
P^\lambda _{\gamma _0,\nu  _1},P^{\prime\bar\lambda }_{\gamma _0,\nu '
_1}:
H^{s-\frac12}(\Sigma )\to H^{s-\frac32}(\Sigma ),\label{tag4.33}
\end{equation}
when $\tfrac32
<s\le 2$. Let us also introduce the operator $\nu _1''=\nu _1'+\mathcal
A'_0\gamma _0$, then Green's formula \eqref{tag2.1b} takes the form
\begin{equation}
(Au,v)_{L_2(\Omega )}-(u,A'v)_{L_2(\Omega )}=(\nu  _1u,\gamma _0v)_{\frac12,-\frac12}-
(\gamma _0u,\nu ''_1v)_{\frac32,-\frac32},\label{tag2.1c}
\end{equation}
for $u\in H^2(\Omega )$, 
$v\in D(\Ama')$, and
 $P^{\prime\bar\lambda }_{\gamma _0,\nu ''
_1}$ (cf.\ \eqref{tag2.8a}) likewise maps as in
\eqref{tag4.33}\,ff. Applying \eqref{tag2.1c} to functions $u,v$ with
$Au=0$, $A'v=0$, we see that 
$P^\lambda _{\gamma _0,\nu  _1}$ and $P^{\prime\bar\lambda }_{\gamma _0,\nu ''
_1}$ are contained in each other's adjoints. 
Therefore
$P^\lambda _{\gamma _0,\nu  _1}$ considered in \eqref{tag4.33} has the
extension $(P^{\prime\bar\lambda }_{\gamma _0,\nu ''_1})^*$, which is
continuous from $H^{s'+\frac32}(\Sigma )$ to $H^{s'+\frac12}(\Sigma )$
for $-2\le s'<-\frac32$. This extends the statement in \eqref{tag4.33} to the
values $0\le s <\frac12$, and by interpolation we obtain the validity
of \eqref{tag4.33} for $0\le s\le 2$. 

$P^\lambda _{\gamma _0,\nu _1}$ can in the localizations to $\rnp$ be
described as the composition of the operator $\nu _1=s_0\gamma _1$ (with
$s_0\in C^{0,1}$) and a generalized Poisson operator consisting
of an operator in order-reduced $x$-form having $C^{0,1}$-smoothness plus a remainder of
lower order. For $s\in \,]\frac32,2]$ we can apply Theorem \ref{Theorem3.1}
$2^\circ$ to the compositions, using that $K^\lambda _\gamma $ is locally the sum
of a composition $\Lambda ^{-2}_{-,+}K_1^0(\lambda )\Lambda _0^2$
(multiplied with smooth cut-off functions) where $K_1^0(\lambda )$
is in $x$-form, and a remainder of lower order.
This implies that $P^\lambda _{\gamma _0,\nu _1}$, apart
from the remainder term coming from $K^\lambda _\gamma $, is the
sum of a first-order
$\psi $do in $x$-form with $C^{0,1}$-smoothness and a remainder term,
mapping $H^{t+1}(\Sigma )$ to $H^t(\Sigma )$ for $|t|<1$, resp.\ 
$H^{t+1-\theta }(\Sigma )$ to $H^t(\Sigma )$ for $-1+\theta
<t<1$. With $s-\frac12=t+1$, $s$ runs in $\,]\frac12, \frac52[\,$
resp.\ $\,]\frac12+\theta , \frac52[\,$ here, which covers the
interval $s\in \,]\frac32,2]$ allowed by the other remainder.  

For low values of $s$ there is again the difficulty that we are
dealing with a composition with ingredients of order-reducing
operators and $x$- or $y$-form operators, which each have 
different rules
for the spaces in which the decompositions and mapping properties 
are valid, and we
refrain from a further discussion here.

Observe moreover that $P^\lambda _{\gamma _0,\nu _1}$ is elliptic (the principal symbol is
invertible) --- since this is known for $P^0_{\gamma _0,\gamma _1}$
(\cite{A62}, \cite{G71}). 

This shows:

\begin{theorem}\label{Theorem4.3} Assumptions as in
  Theorem {\rm \ref{Theorem4.2}}. $P^\lambda _{\gamma _0,\nu  _1}$ and $P^{\prime\bar\lambda }_{\gamma _0,\nu '
_1}$ 
map
$H^{s-\frac12}(\Sigma )$ to $H^{s-\frac32}(\Sigma )$ for $s\in [0,2]$.
They are generalized elliptic $\psi $do's of order $1$, in the sense
that for $s\in \,]\frac32,2]$, they have the form of an elliptic principal part in $x$-form of H\"older smoothness $(0,1)$ plus a
lower order part.
%; for $s\in [0,\frac12[\,$, they can be written as
%the sum of a principal part in $y$-form of H\"older smoothness $(0,1)$ and a
%lower order part. 
\end{theorem}

With these mapping properties it is straightforward to verify 
that $\Gamma $ and $\Gamma '$
defined in \eqref{tag2.9a} satisfy the full statement in \eqref{tag2.9}.

When more smoothness of $\Omega $ is
assumed, the representation of $P^\lambda
  _{\gamma _0,\nu _1}$ as the sum of a principal part in $x$-form and
  a lower-order term can of course be extended to larger intervals
  than found above.

\section{Interpretation of realisations.}\label{Section5}

We now have all the ingredients to interpret the abstract
characterisation of closed realisations $\wA$ in terms of operators
$T:V\to W$ recalled in Section \ref{Section1}, to boundary conditions. In fact,
we have the mappings defined from the trace operator $\gamma _0$ 
\[
\gamma _{Z_\lambda }: Z_\lambda \simto H^{-\frac12}(\Sigma ),\quad \gamma _{Z_\lambda }^*:
H^\frac12(\Sigma )\simto Z_\lambda ,
\]
and the  mappings defined from Poisson-type operators
\[
K^\lambda _\gamma :H^{-\frac12}(\Sigma )\to H^0(\Omega ),\quad
(K^\lambda _\gamma )^*:H^0(\Omega )\to H^{\frac12}(\Sigma ),
\]
as well as the versions with primes. Then the various definitions
recalled in Section \ref{Section2} for the smooth case, carrying $T^\lambda :V_\lambda \to
W_{\bar\lambda }$ 
over to $L^\lambda :H^{-\frac12}(\Sigma )\to H^\frac12(\Sigma )$ if
$V=Z$, $W=Z'$, resp.\ to $L_1^\lambda :X_1\to Y_1$ in general, are
effective in exactly the same way, and all the diagrams are valid in this situation.

In this way, $\wA$ is determined by a Neumann-type boundary condition
\[
\nu _1u=(L+P^0_{\gamma _0,\nu _1})\gamma _0u
\]
in the case $V=Z$, $ W=Z'$, and by a condition involving projections
 in the general case.

The adjoint $\wA$ is determined by the boundary condition 
\[
\nu '_1u=(L^*+P^{\prime 0}_{\gamma _0,\nu '_1})\gamma _0u
\]
in the case $V=Z$, $ W=Z'$ (resp.\ by a condition involving projections
in the general case), where $L^*$ is the adjoint of $L$, considered as
a generally unbounded operator from $H^{-\frac12}(\Sigma )$ to
$H^{\frac12}(\Sigma )$. 

There is a well-defined $M$-function $M_L(\lambda )$, which coincides with
$-(L^\lambda )^{-1}$ for $\lambda \in \varrho (A_\gamma )\cap \varrho
(\wA)$; here \eqref{tag2.18a} and \eqref{tag2.19} hold. Suitably
modified results hold in cases of general $V,W$.

For the case $V=Z$, $W=Z'$, we have obtained:

\begin{theorem}\label{Theorem5.1} When $\Omega $ is $C^{1,1}$ and
  $A$ has $C^\infty $ coefficients, bounded with bounded derivatives on
  a neighborhood of $\Omega $, and is uniformly strongly elliptic,
  then Theorem {\rm \ref{Theorem2.1} (i)--(v)} and {\rm \eqref{tag2.18a}}
  are valid.

\end{theorem}

Gesztesy and Mitrea have in \cite{GM08} established Kre\u\i{}n resolvent
formulas for the Laplacian under a weaker smoothness hypothesis, namely
that $\Omega $ is $C^{1,\sigma }$ with $\sigma  >\frac12$. Here they
treat {\it selfadjoint} realisations determined by Robin-type boundary
conditions 
\begin{equation} \gamma _1u=B\gamma _0u,\label{tag5.1}
\end{equation}
with $B$ compact from $H^1$ to $H^{0}$ (assured if $B$ is of order
$<1$). Posilicano
and Raimondi \cite{PR08} describe results for {\it selfadjoint} realisations
in case $\Omega $ is $C^{1,1}$ and the  coefficients of
$A$, when it is written in symmetric divergence form, are $C^{0,1}$ satisfying
various hypotheses. They remark that their treatment works for
boundary conditions  \eqref{tag5.1} with $\gamma _1$ replaced by the
oblique Neumann trace operator $\nu _A$ \eqref{tag2.1e} connected with the 
divergence form. Here
$B$ is taken of order $<1$ and elliptic (we do not quite see the relevance
of the latter hypothesis), so it is a Robin-type perturbation of the natural Neumann
condition.

It is an important point in the present treatment, besides that it
deals with nonselfadjoint situations, that Neumann-type conditions
\eqref{tag2.12} with general $\psi $do's $C$ of order 1 are
included in the detailed discussion.

Furthermore, our pseudodifferential strategy  allows the application of ellipticity
concepts:

When $C$ is a generalized pseudodifferential operator of order 1 and H\"older
smoothness $(0,1)$, $L=C+P^0_{\gamma _0,\nu _1}$ is a generalized
pseudodifferential operator of order 1 and H\"older smoothness $(0,1)$, and
vice versa. $L$
is elliptic precisely when the model boundary value problem for $A$
with the boundary condition \eqref{tag2.12} is uniquely solvable at
all $(x',\xi ')$ with $\xi '\ne 0$ in the boundary cotangent space
(this is the Shapiro-Lopatinski\u\i{} condition). $L^\lambda $ is then
also elliptic at each $\lambda \in \varrho (A_\gamma )$ (since
$P^\lambda _{\gamma _0,\nu _1}-P^0_{\gamma _0,\nu _1}$ is of order
$<1$). 

Moreover, there is then a parametrix of
$L$, and this can be used to investigate the regularity of the domain
of $L$.
Likewise, each $L^\lambda $ has a parametrix
then. However, we 
 want to set the  true inverse $-M_L(\lambda )$ in relation to such a
parametrix. 

Restrict the attention to the case where
$C$ is a first-order {\it differential} operator on $\Sigma $ with $
C^{0,1}$-coefficients; then we can say more about $M_L(\lambda )$ with the
present methods.

Assume a little more, namely that there is a ray $\lambda =-\mu ^2e^{i\theta
}$, $\mu \in{\Bbb R}$, such that when we include $\lambda $ in the 
principal symbol of 
$P^\lambda _{\gamma _0,\nu _1}$, then the principal symbol of
$L^\lambda =C-P^\lambda _{\gamma _0,\nu _1}$ is invertible for $|\xi '|^2+|\mu
|^2\ge 1$ (``parameter-ellipticity''). 
Let $s\in \,]\frac32,2]$. As in Section \ref{Section4},
we can invoke the system for $\widehat A$ on $\widehat\Omega =\Omega \times S^1$
\eqref{4.15b} coupled with the same boundary operator (constant
in the $t$-direction)
\begin{equation}
\widehat {\mathcal A} =\begin{pmatrix} \widehat A\\ \nu_1-C\gamma
_0\end{pmatrix}:H^{s}(\widehat\Omega ) \to\begin{matrix}
H^{s-2}(\widehat\Omega )\\ \times \\ H^{s-\frac32}(\widehat\Sigma )
\end{matrix}; \label{tag5.2}
\end{equation}
it is elliptic and has a parametrix $\widehat{\mathcal B}^0 $. For
the functions $u(x,t)=w(x)e^{i\mu t}$, this gives a $\lambda
$-dependent parametrix family for $\mathcal A(\lambda )=\begin{pmatrix}  A-\lambda \\ \nu_1-C\gamma
_0\end{pmatrix}$ (when \linebreak$|\lambda |\ge 1$) such that
the remainder in the composition with 
$\mathcal A(\lambda )$ is $O(\ang\mu ^{-\theta })$ for $\lambda \to \infty $ on
the ray. Then there is a true inverse of $\mathcal A(\lambda )$, hence of  $L^\lambda $, for sufficiently large $\lambda
$ on the ray. We can follow this up for the operator $\widehat
L=C-\widehat P_{\gamma _0,\nu _1}$ over $\widehat\Sigma $, which gives
$L^\lambda $ when applied to functions $\varphi (x')e^{i\mu t}$. Here
$\widehat L$ has a parametrix $\widehat{\widetilde L}$ such that
$\widehat L\widehat{\widetilde L}-I$ is of negative order; this gives
a parametrix $\widetilde L^\lambda $ of $L^\lambda $ such that
$ L^\lambda \widetilde L^\lambda -I$ has an $O(\ang\mu ^{-\theta })$
estimate.
For sufficiently large $\lambda $ on the ray this allows us to write   
$M_L(\lambda )=-(L^\lambda )^{-1}$ as $-\widetilde L^\lambda +\mathcal
R$ with $\mathcal R$ of lower order. More precisely, $\widetilde
L^\lambda $ is obtained as a composition of an operator in $x$-form
with an order-reducing operator to the left; it maps
from
$H^{s-\frac32}$ to $H^{s-\frac12}$, and the remainder maps
from $H^{s-\frac32-\theta }$ to $H^{s-\frac12}$.  (The $s\in
\,]\frac32,2]$ run inside the interval where the parametrix
construction for elliptic first-order $\psi $do's of H\"older
smoothness $(0,1)$ works, as in Theorem \ref{Theorem3.1} $3^\circ$ and Remark \ref{Remark3.2}.) In this sense, $M_L(\lambda )$ is a
generalized $\psi $do of order $-1$.

Using this information for $s=2$, we see that 
$M_L(\lambda )$ map $H^\frac12$ not
just to $H^{-\frac12}$, but to $H^\frac32$. Then $D(L)=D(L^\lambda
)=H^{\frac32}$ and $D(\wA)$ is in $H^2(\Omega
)$. 

If, moreover, $C^*$ has H\"older smoothness $C^{0,1}$,  
the adjoint $\wA^*$ is of the same type. In particular, there is 
selfadjointness if $A$ and $L$ are formally selfadjoint. This gives a
very satisfactory version of the Kre\u\i{}n formula.

\begin{theorem} \label{Theorem 5.2} If, in addition to the hypotheses
  of Theorem {\rm \ref{Theorem5.1}}, $C$ is a first-order differential
  operator with H\"o{}lder smoothness $(0,1)$ and the principal symbol of $L^\lambda =C-P^\lambda _{\gamma _0,\nu
    _1}$ is parameter-elliptic on a ray $\lambda =-\mu ^2e^{i\theta
}$, $\mu \in{\Bbb R}$, then $D(L)=H^{\frac32}(\Sigma )$, and
$M_L(\lambda )$ is for large $\lambda $ on the ray the sum of  an
elliptic  $\psi $do of order $-1$ and H\"older smoothness $(0,1)$, in
order-reduced $x$-form, and a lower-order term.
Then $D(\wA)\subset H^2(\Omega
)$.

If, moreover, $C^*$ has H\"older smoothness $(0,1)$, 
the adjoint $\wA^*$ is defined similarly from of $L^*$ with
$D(L^*)=H^{\frac32}$,  $D(\wA^*)\subset H^2(\Omega
)$. In particular, $\wA$ is
selfadjoint if $A$ and $L$ are formally selfadjoint. 
\end{theorem}

From the point of view of the systematic parameter-dependent calculus
of \cite{G96}, the
 symbols of $C$
and $P^\lambda _{\gamma _0,\nu _1}$ have ``regularity $\nu = +\infty
$'' when $C$ is a differential operator, so there is a parametrix with the same
``regularity $+\infty $''.

{\it Pseudodifferential} operators $C$ can be included in the discussion if
the symbol classes in \cite{G96} are used in a more definitive way
(here when $C$ is of order 1, it has ``regularity 1'', and the same
will hold for the resulting principal symbols of $L^\lambda $ and
$M_L(\lambda )$). Considerations with finite positive ``regularity''
play an important role in \cite{A03, A05b}. We hope to return to such
cases in future works, but here just wanted to show what can be 
done using Agmon's
principle.

\bigskip

\begin{flushleft}

{\bf AMS Subject Classification: 35J25, 47A10, 58J40.}\\[2ex]

% Write more than one author separately if they have different 
% affiliations, otherwise write the names on the same line, separeted 
% by commas.
%
Gerd GRUBB,
\\ Department of Mathematical Sciences,
\\ Copenhagen University,
\\Universitetsparken 5,
\\2100 Copenhagen, DENMARK
 \\ \texttt{grubb@math.ku.dk}

\end{flushleft}

\end{document}